\documentclass[10pt, oneside]{article}
\usepackage{geometry}
\usepackage{graphicx}
\usepackage{amssymb}
\usepackage{mathtools}
\usepackage{amsrefs}
\usepackage{epstopdf}
\usepackage{enumitem}
\usepackage[utf8]{inputenc}
\usepackage{subcaption}
\usepackage{hyperref}
\usepackage{xspace}
\usepackage{amsthm}

\theoremstyle{remark}
\newtheorem{remark}{Remark}

\newcommand{\Om}{\Omega}

\newcommand*{\eg}{e.g.\xspace}
\newcommand{\dumux}{\mbox{DuMu$^\mathsf{x}$\xspace}}

\newcommand{\e}{\ensuremath{\epsilon}}

\newcommand{\Oe}[1]{\mathcal{O}(\e^{#1})}

\newcommand{\p}{\partial}

\newcommand{\Rey}{\mathrm{Re}\,}
\newcommand{\Ca}{\mathrm{Ca}\,}
\newcommand{\ve}[1]{\ensuremath{\mathbf{#1}}}
\newcommand{\vv}{\ve{v}}
\newcommand{\vw}{\ve{w}}
\newcommand{\vg}{\ve{g}}
\newcommand{\Eu}{\ensuremath{\mathrm{Eu}}}
\newcommand{\Fr}{\ensuremath{\mathrm{Fr}}}

\newcommand{\Po}{\ensuremath{\mathcal{P}}}
\newcommand{\dd}{\,\mathrm{d}}
\newcommand{\x}{\ve{x}}
\newcommand{\y}{\ve{y}}

\newcommand{\porosity}{\ensuremath{\varphi}}
\newcommand{\pf}{\ensuremath{u}}
\newcommand{\mb}{\ensuremath{S}}
\newcommand{\n}{\ensuremath{\ve{n}}}
\newcommand{\Gr}{\ensuremath{\mathcal{G}}}
\newcommand{\al}{\alpha}
\newcommand{\InterfaceFluids}{{\Gamma_f}}
\newcommand{\InterfacePhases}{\Gamma_\al}
\newcommand{\InterfaceCell}{\Gamma}
\newcommand{\contactAngle}{\ensuremath{\theta_\mathrm{eq}}}
\newcommand{\pfBCcell}{- \xi^{-1} \cos(\contactAngle) \sqrt{2P(\pf_0)}}
\newcommand{\slip}{\ensuremath{\lambda}}

\newcommand{\effPerm}{\ensuremath{\mathcal{K}}}
\newcommand{\effSurf}{\ensuremath{\mathcal{M}}}
\newcommand{\ratioVisc}{\ensuremath{M}}
\newcommand{\ratioDens}{\ensuremath{R}}
\newcommand{\dimension}{\mathrm{d}}
\newcommand{\interface}{\ensuremath{\delta}}
\newcommand{\tension}{\ensuremath{\gamma}}
\newcommand{\tens}[1]{\ve{#1}}
\newcommand{\mob}{\ensuremath{\sigma}}

\newcommand{\xxi}{{\bar{\xi}}}
\newcommand{\refL}{L}
\newcommand{\out}{\mathrm{out}}
\newcommand{\inn}{\mathrm{in}}
\newcommand{\s}{\ve{s}}
\newcommand{\Ox}[1]{\mathcal{O}(\xxi^{#1})}
\newcommand{\Interface}{\Gamma}

\begin{document}

\title{Upscaling and Effective Behavior for Two-Phase Porous-Medium Flow using a Diffuse Interface Model}

\author{Mathis Kelm
    \footnote{Institute for Modelling Hydraulic and Environmental Systems, University of Stuttgart, Pfaffenwaldring 61, 70569 Stuttgart, Germany},
    Carina Bringedal
    \footnote{Department of Computer science, Electrical engineering and Mathematical sciences, Western Norway University of Applied Sciences, Inndalsveien 28, 5063 Bergen, Norway}
    \ \footnotemark[1],
    Bernd Flemisch\footnotemark[1]
}

\maketitle

\begin{abstract}
We investigate two-phase flow in porous media and derive a two-scale model, which incorporates
pore-scale phase distribution and surface tension into the effective behavior at the larger Darcy
scale.
The free-boundary problem at the pore scale is modeled using a diffuse interface approach in the
form of a coupled Allen-Cahn Navier-Stokes system with an additional momentum flux due to surface
tension forces. Using periodic homogenization and formal asymptotic expansions, a two-scale model
with cell problems for phase evolution and velocity contributions is derived. We investigate the
computed effective parameters and their relation to the saturation for different fluid
distributions, in comparison to commonly used relative permeability saturation curves.
The two-scale model yields non-monotone relations for relative permeability and saturation.
The strong dependence on local fluid distribution and effects captured by the cell problems
highlights the importance of incorporating pore-scale information into the macro-scale equations.
\end{abstract}

\section{Introduction}\label{sec1}

Flow through porous media, especially in multi-phase systems, is of interest in a variety of
applications from oil recovery and $CO_2$ sequestration to fuel cells and biological systems.
Traditionally these are modeled using extended Darcy's law for two phases, which lacks the
mathematical derivation from pore-scale information available for the single-phase Darcy's law.
Single-phase Darcy's law can be derived from a pore-scale model of
(Navier-)Stokes equations through upscaling procedures. For two phases Darcy's law has been extended
by introducing empirically derived relative permeability saturation curves to capture fluid interactions.
The effective behavior depends only on the averaged saturation and the model fails to capture
different pore-scale effects.
Further empiric modifications such as play-type hysteresis have been added to account for missing
behavior.
In this work we start from a two-phase-flow model at the pore scale and through appropriate assumptions and a
periodic homgenization approach, we derive a two-scale model of macroscopic Darcy's law-type equations with effective
parameters computed from solutions of resulting pore-scale cell problems.

The flow of the two fluids is modeled on the pore scale using quasi-compressible Navier-Stokes equations with phase-dependent
density and viscosity as well as an additional term
capturing surface tension forces at the fluid-fluid interface.
For the pore-scale model with resolved phase distribution we use a diffuse interface approach in the form
of an Allen-Cahn phase-field model \cite{allen1979}.
Where a sharp-interface model separates the pore space into subdomains for each fluid phase, our
approach captures the interface implicitly through a smooth phase field function.
All unknowns are defined and modified equations hold on the entire domain of the pore space.
This removes the need to track or capture the interface and
decompose the domain in the numeric simulation.
To resolve the transition zone between the two bulk phases a fine discretization is required near
the interface, requiring high computational effort in the absence of adaptive local mesh refinement.
Furthermore the Allen-Cahn model used in this work is not conservative, including a curvature driven
motion of the interface.
While the conservative Cahn-Hilliard model requires solving a fourth order partial differential
equation, we consider the second order Allen-Cahn equation, which offers a much simpler numerical
implementation, and we investigate the two-scale model derived from it.
While there are different approaches to ensuring conservation of the phase field
\cite{bringedal2020con} or counteracting this curvature-driven motion entirely \cite{xu2012}, we
apply the scaled mobility approach presented in \cite{abels2021non} to eliminate the
curvature-driven motion only in the sharp-interface limit.
A core advantage of the diffuse interface approach is the ease of upscaling the equations, as they
are defined on
a stationary domain encompassing both fluid phases.
In contrast, upscaling a sharp-interface model requires special attention to the evolving domains
\cite{vanNoorden2009}.

The model is presented for two phases but extends well to more evolving phases
\cites{kelm2022, redecker2016}, only a non-neutral contact angle at internal triple points between
phases requires more attention \cite{rohde2021}.

To bridge the scale gap between pore scale and the averaged Darcy scale we employ the method of formal asymptotic homogenization
\cite{hornung1997homogenization}.
Representing the porous medium by a periodic arrangement of scaled reference cells, one introduces
an additional spatial coordinate for this smaller scale.
Here the macroscopic scale corresponds to the Darcy scale and the microscopic scale is the pore
scale.
These two scales are sufficiently separated, with representative length scales $L$ and $\ell$
respectively defining a length scale ratio $\e = \ell/L \ll 1$.
Assuming asymptotic expansions of the unknowns in terms of $\e$, considering the
limit $\e \to 0$ and grouping by orders of $\e$, one obtains equations that are defined on the macro
scale and form micro-scale cell problems. The former contains effective parameters which are computed
from the cell problem unknowns, linking the two scales.
Through the cell problems, these parameters and the effective behavior depend
additionally on the distribution of
fluids and are able to capture more detailed effects on the fluid velocities.

The homogenization leads to macroscopic equations similar to the extended Darcy's law but with
effective parameters computed from solutions of cell problems instead of being given by empiric
relations.
These effective parameters can be understood as total phase mobilities and for isotropic pore
geometries a relative permeability can be computed.
Through the cell problems these relative permeabilities in the two-scale model depend on the local
phase distribution at the pore scale.
This motivates a comparison of computed relative permeabilities to commonly used functional relations,
enabling us to investigate under which conditions the models agree and which additional pore-scale effects are
captured by the two-scale model.
We investigate numerically the effective parameters computed for different geometries and fluid
properties. Constructing fluid distributions of varying saturations, we solve the cell problems for
the velocity contributions to obtain a relative permeability curve.
The models are implemented using a staggered finite volume discretization in
\dumux \cite{KOCH2021423}, with
Dune-SPGrid \cite{nolte2011} and Dune-Subgrid \cite{graeser2009} used for the grid.

In \cite{daly2015}, \cite{metzger2021}, and \cite{sharmin2022} a Cahn-Hilliard model was used to
derive a similar two-scale model. We instead
investigate the applicability of the simpler Allen-Cahn model and the effects of micro-scale
information on the effective flow parameters.
The derivation of a two-scale model in \cite{metzger2021} additionally assumes a separation of time
scales and a different asymptotic expansion to arrive at an instationary problem for the phase
distribution. We use the full expansion in the spatial separation \(\e\) instead and afterwards introduce an
artifical evolution to obtain a local phase distribution.
In \cite{daly2015} three time scales are used to separate interface equilibration, macroscopic flow
and saturation changes.
In \cite{sharmin2022} a solute-dependent surface tension is considered and simulations of the
coupled two-scale model are presented. We instead focus on a systematic investigation of the
effective parameters computed from the cell problems. 
A similar comparison of effective parameters is included in \cite{lasseux2022}, where a
sharp-interface Stokes model for two-phase flow is upscaled using volume averaging.
However, \cite{lasseux2022} does not account for a three-phase contact line or slip velocities at
the solid, which are accounted for in the current work.

In this contribution we aim to derive a two-scale model for two-phase flow in porous media, which
captures effects of pore-scale fluid distribution and surface tension on the effective behavior at
the macro-scale. Our goal is to determine conditions under which the computed effective parameters
show significantly different behavior to commonly used saturation-dependent curves and highlight the
importance of resolving pore-scale behavior through a two-scale model for two-phase porous-medium
flow.

The structure of this paper is as follows.
We present our pore-scale model for two-phase flow using an Allen-Cahn formulation
coupled to a modified Navier-Stokes system in Section~\ref{sec:pore_scale}.
In Section~\ref{sec:upscaling} we perform the upscaling and obtain
a two-scale model using periodic homogenization.
We investigate the
dependence of computed effective parameters on pore-scale conditions using \dumux in Section~\ref{sec:numerics}.
Section~\ref{sec:conclusion} summarizes the derived model and relates it to the extended Darcy's
law.

\section{Pore-scale modeling}\label{sec:pore_scale}
We begin by presenting the pore-scale model describing two-phase flow. The model
is first presented in a sharp-interface formulation before we introduce the diffuse-interface
description for the phase-field model.

The general domain is a stationary void space inside a porous medium.
We consider a domain $\Om = \Po \cup \Gr$ decomposed into a domain $\Po$ capturing the pore space
and the solid matrix $\Gr$.
Here the distribution of two fluid phases needs to be captured and evolved according to the flow of
the fluids, including effects of surface tension at the fluid-fluid interface.

\subsection{Sharp-interface model}
\label{sec:si_model}
\begin{figure}
\center
\includegraphics[width=0.5\textwidth]{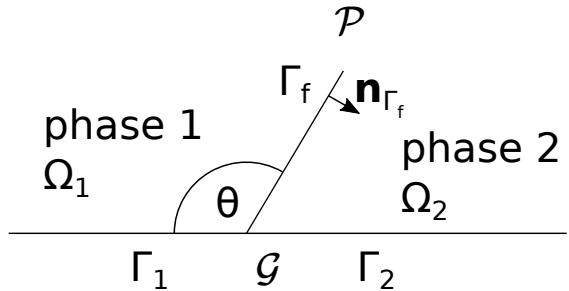}
\caption{Interfaces and triple point.}
\label{fig:domain_sharp}
\end{figure}
The model is first formulated with a sharp interface separating the fluid phases and corresponding
evolving subdomains $\Om_\al(t) \subset \Po$, $\al = 1,2$ of the stationary pore space $\Po$
(Figure~\ref{fig:domain_sharp}).
The interfaces between the fluid phases, denoted by $\InterfaceFluids(t) = \overline{\Om_1(t)} \cap
\overline{\Om_2(t)}$, as well as the fluid-solid interfaces
$\InterfacePhases = \overline{\Om_\al} \cap \Gr$ evolve accordingly.
Within each subdomain the flow is modelled using incompressible Navier-Stokes equations, with
constant viscosity $\mu_\al$ and density $\rho_\al$.
\begin{align}
\label{eq:si_mass}
\nabla\cdot \vv_\al &= 0
~, && \text{ in }\Om_\al ~,\\
\label{eq:si_momentum}
\rho_\al \frac{\p}{\p t} \vv_\al + \rho_\al \nabla\cdot (\vv_\al \otimes \vv_\al)
&= -p_\al + \mu_\al \nabla^2 \vv_\al + \rho_\al \ve{g}
~, && \text{ in }\Om_\al ~.
\end{align}
At the interfaces $\InterfacePhases$ with the solid matrix we prescribe a Navier-Slip
\cite{thalakkottor2016}
boundary condition with slip length $\slip \ge 0$ for the velocity $\vv_\al$
\begin{equation}
\label{eq:si_bc_slip}
\vv_\al = - \slip ( \partial_\n \vv_\ve{t} ) \ve{t}
~,
\qquad
\text{on }\InterfacePhases ~,
\end{equation}
where $\ve{t}$ is tangential to the fluid-solid interface $\partial\Gr$.
At the fluid-fluid interface $\InterfaceFluids$ we prescribe continuity of velocities with a no-slip
condition enforcing zero tangential velocity, with interface velocity $V_\InterfaceFluids$ in the
direction of the interface normal $\n_\InterfaceFluids$ pointing from fluid 1 into fluid 2,
\begin{align}
\vv_\al = V_\InterfaceFluids \n_\InterfaceFluids ~.
\end{align}
The interface itself moves due to the local velocity of the fluids
and with surface tension $\tension$ and curvature $H_\InterfaceFluids$
we have a jump condition
\begin{equation}
\label{eq:si_interface_jump}
[-p_\al \tens{I}
+ \mu_\al (\nabla \vv_\al + (\nabla \vv_\al)^T - \frac{2}{3} \nabla \cdot \vv_\al \tens{I})
]\n_\InterfaceFluids
= - \tension H_\InterfaceFluids \n_\InterfaceFluids
\qquad \text{on }\InterfaceFluids ~,
\end{equation}
where $[\psi] = \psi_1 - \psi_2$ denotes the jump over the interface.
At the three-phase contact point the fluid-fluid-interface meets the solid with a contact angle
$\contactAngle$ (Figure~\ref{fig:domain_sharp}).
In a sharp-interface formulation this could be incorporated using a boundary condition for a level
set equation, when a level set is used to track the location of the fluid-fluid-interface.

\subsection{Phase-field model}
\label{sec:pf_model}
To model two-phase flow at the pore scale we use a diffuse-interface approach, capturing the
distribution of phases with a phase-field function $\pf$ defined on the total void space
$\pf : T \times \Po \to [0, 1]$ with $\pf=0$ and $\pf=1$ corresponding to the two distinct
phases $\al=2$ and $\al=1$ respectively.
The phase-field variable evolves according to an advective Allen-Cahn equation, a second order
partial differential equation (PDE) with a non-linear source term.
The multi-phase system is then modeled as one fluid with varying properties, depending on
$\pf(t, \x)$ to account for phase distribution.

We use the compressible Navier-Stokes equations for
non-constant viscosity, derived from conservation laws using Stokes hypothesis $\zeta=0$.
A surface tension flux is added to the momentum equation, coupling it to the Allen-Cahn equation as
presented for the stationary Stokes equation in \cite{abels2018}. In \cite{abels2012} the surface
tension term is introduced to a Navier-Stokes equation with a phase field evolved according to the
Cahn-Hilliard equation.

The resulting Navier-Stokes equations in conservative form are
\begin{subequations}
\begin{align}
\label{eq:diff_mass}
\frac{\partial \rho}{\partial t} + \nabla \cdot (\rho \vv) &= 0
&\text{ in }\Po~,\\
\label{eq:diff_momentum}
\frac{\partial}{\partial t}(\rho \vv) + \nabla \cdot (\rho \vv \otimes \vv) &=
 -\nabla p
 + \nabla \cdot \left( \mu \left(\nabla\vv + (\nabla\vv)^T
             - \frac{2}{3} (\nabla\cdot\vv) \mathbf{I} \right) \right)
         + \rho \vg
& \\ &\phantom{=} \notag
         - \frac{3}{2} \xi \tension \nabla\cdot (\nabla\pf \otimes \nabla\pf)
&\text{ in }\Po~.
\end{align}
\end{subequations}
We focus on the case of each fluid being incompressible, treating the multiphase system as a
quasi-incompressible fluid with density and viscosity dependent on the present phase in a linear
manner.
Denoting mass and viscosity ratios $R = \rho_1 / \rho_2$ and $M = \mu_1 / \mu_2$, we write
\begin{subequations}
\label{eq:fluid_properties_linear}
\begin{align}
\rho(\pf) &= \rho_1 \pf + \rho_2 (1-\pf) = \rho_2 + \pf (\rho_1 - \rho_2)
    = \rho_2 \left(1+\pf\left(\frac{\rho_1}{\rho_2} - 1\right)\right) = \rho_2 (1 + \pf (R-1) )
~, \\
\mu(\pf) &= \mu_1 \pf + \mu_2 (1-\pf) = \mu_2 ( 1 + \pf (M-1) )
~.
\end{align}
\end{subequations}
This is combined with the advective Allen-Cahn equation, additionally coupled via the advecting velocity and the
surface tension momentum flux.
The phase-field equation is given as
\begin{equation}
\label{eq:pore_phasefield}
\frac{\p \pf}{\p t} + \nabla \cdot (\vv \pf) = \mob \xi ( \nabla^2 \pf - \xi^{-2} P'(\pf) )
\,,\quad
\text{in }\Po
\,,\quad
P(\pf) = 8 \pf^2 (1-\pf)^2
~,
\end{equation}
where we consider the scaled diffusivity $\mob \xi$ suggested in \cite{abels2021non} in order to
obtain a more favorable sharp-interface limit.
The diffusivity $\mob$ and diffuse-interface width $\xi$ prescribe how strongly a certain profile
is enforced for the transition zone between bulk phases and the width of this zone respectively.
The double-well potential, $P(\pf)$, encourages a separation of phases through stable minima at $\pf = 0$
and $\pf = 1$.

\begin{figure}
\center
\includegraphics[width=0.5\textwidth]{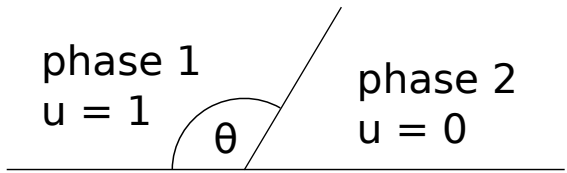}
\caption{Contact angle boundary condition for the diffuse-interface phase-field model.}
\label{fig:contact_angle_diff}
\end{figure}

The (no-)slip boundary conditions on solid boundaries apply also here, with interpolation based on
the phase-field $\pf$ if phase-dependent slip lengths are desired.
Additionally for the phase-field a Neumann boundary condition encoding the contact angle
$\contactAngle$ measured through fluid 1, corresponding to $\pf = 1$
(Fig~\ref{fig:contact_angle_diff}), is
prescribed \cites{frank2018, lunowa2021}.

\begin{subequations}
\begin{alignat}{2}
\label{eq:diff_bc_slip}
\vv &= - \slip (\partial_\n \vv_\ve{t}) \ve{t}
~, && \text{ on }\partial\Po
~,\\
\label{eq:diff_bc_contact_angle}
\partial_\n \pf &= -\cos(\contactAngle) \xi^{-1} \sqrt{2P(\pf)}
~, && \text{ on }\partial\Po
~.
\end{alignat}
\end{subequations}
Note that for vanishing slip length $\slip =0$ and for neutral contact angle
$\contactAngle = \frac{\pi}{2}$, the boundary conditions reduce to homogeneous Dirichlet and Neumann conditions
for velocity and phase field respectively.
The Appendix~\ref{app:si} includes the sharp-interface limit of the coupled Navier-Stokes
Allen-Cahn system presented in this section, yielding again the sharp-interface model presented in
Section~\ref{sec:si_model}.

\section{Upscaling using periodic homogenization}\label{sec:upscaling}

Starting from the diffuse-interface model from Section~\ref{sec:pf_model}, describing the phase
distribution and flow at the pore scale, we
derive effective equations capturing the behavior at the Darcy scale. At the same time we obtain
a set of micro-scale cell problems, defined at each macroscopic point, the solutions of which allow
the computation of effective tensors. These are then used as parameters to the
macro-scale problem, which in turn supplies global pressure gradients and saturations used by the
cell problems, arriving at a coupled two-scale problem.
As will be shown, these effective parameters can be interpreted as phase mobilities.

\begin{figure}
\center
\includegraphics[width=\textwidth]{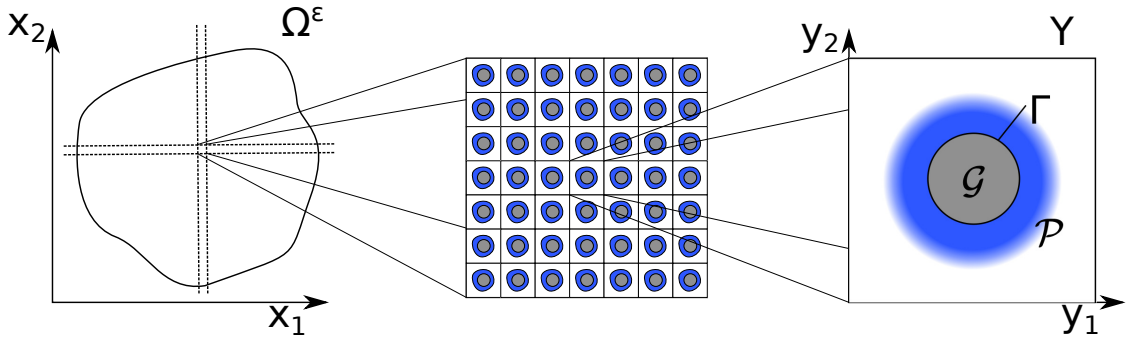}
\caption{Separation of scales and domain definitions}
\label{fig:homogenization}
\end{figure}
To achieve this we assume a sufficient separation of scales, with characteristic length scales $\ell$
for the pore scale and $L$ for the Darcy scale at which we are interested in the effective behavior.
Assuming a scale separation $\e = \ell/L \ll 1$, we first non-dimensionalize the model.
We then follow the approach of periodic homogenization, modeling the porous medium as a
periodic arrangement of scaled reference cells $\e Y$, $Y = [0, 1]^\dimension$ with dimension
$\dimension$, where $\dimension$ is 2 or 3, (see Figure~\ref{fig:homogenization}).
The domain of the porous medium is thus decomposed into
$\Om^\e = \cup_{w \in W_\Om} \,\e(w + Y)$, with
$W_\Om \in \mathbb{Z}^\dimension$ a set of indices.
With local pore space and solid matrix denoted as
$Y = \Po \cup \Gr$ and the outer boundary of the reference cell given as $\partial Y$
we denote the interior boundary between solid and fluid as
\begin{equation}
\InterfaceCell = \overline{\Po} \cap \Gr
~.
\end{equation}
The global pore space and its boundary with the solid matrix in the non-dimensional case is then given as
\begin{equation}
\Po^\e = \cup_{w \in W_\Om} \,\e(w + \Po)
    ~,\quad
\InterfaceCell^\e = \cup_{w \in W_\Om} \,\e(w + \InterfaceCell)
    ~.
\end{equation}
While their fluid content can generally vary between the cells, we here
assume that the porous matrix is constant both in time and space in order to facilitate the
derivation of the two-scale model.
\newcommand{\poro}{\varphi}
As the reference cell is the unit cube, we can define the constant porosity
$\poro = |\Po|$.
The conditions inside the pore spaces still vary, with fluid
saturations and distributions changing in time and between cells.

We use the scale separation to introduce a micro-scale coordinate, rewrite spatial derivatives
accordingly, assign scalings in terms of $\e$ for the dimensionless numbers and assume all unknowns
to have an asymptotic expansion in $\e$.
Inserting the expansions into the model equations and gathering terms of equal order in terms of $\e$ we
obtain a new set of equations containing 
derivatives with respect to coordinates of only one of the spatial scales.
This procedure yields macroscopic equations capturing the effective behavior and cell problems defined on the
reference cell.
Effective parameters are defined through cell problems and integrals of their solutions, linking the two sets of
equations and capturing detailed effects of local phase distribution on effective behavior.
In addition to pressure driven flow, a separate cell problem captures flow due to surface tension
forces, reflected in an added velocity contribution on the effective scale.

\subsection{Non-dimensionalization}
\label{sec:non-dim}
In preparation of upscaling by periodic homogenization we non-dimensionalize the micro-scale model.
The homogenization requires a separation of scales between representative length $\ell$ at the pore scale
and the length scale of interest $L$ at the macro scale, quantified by the small number
$\e = \ell/L \ll 1$.
Defining reference values with dimensions
\begin{align*}
[\hat{L}] &= \mathrm{m}
&
[\hat{\ell}] &= \mathrm{m}
&
[\hat{\xi}] &= \mathrm{m}
&
[\hat{t}] &= \mathrm{s}
&
[\hat{v}] &= \mathrm{\frac{m}{s}}
\\
[\hat\rho] &= \mathrm{\frac{kg}{m^3}}
&
[\hat\mu] &= \mathrm{\frac{kg}{m \cdot s}}
&
[\hat{p}] &= \mathrm{\frac{kg}{m \cdot s^2}}
&
[\hat\slip] &= \mathrm{m}
&
[\hat\mob] &= \mathrm{\frac{m}{s}}
~,
\end{align*}
and letting
\begin{align}
\label{eq:non-dim_reference}
\hat{L} &= L
&
\hat{\ell} &= \ell
&
\hat\xi &= \hat{\ell}
&
\hat{t} &= \frac{\hat{L}}{\hat{v}}
&
\hat\rho &= \rho_2
&
\hat\mu &= \mu_2
&
\hat\slip &= \hat{\ell}
&
\hat\mob &= \mob
~,
\end{align}
one can rewrite the phase-field model equations using non-dimensional variables and parameters
\begin{align*}
\bar\xi &= \frac{\xi}{\hat{\xi}}
&
\bar\vv &= \frac{\vv}{\hat{v}}
&
\bar\rho &= \frac{\rho}{\hat\rho}
&
\bar\mu &= \frac{\mu}{\hat\mu}
&
\bar{p} &= \frac{p}{\hat{p}}
&
\bar{t} &= \frac{t}{\hat{t}}
&
\bar\slip &= \frac{\slip}{\hat{\slip}}
~.
\end{align*}
In the following all variables and parameters are non-dimensional and the overline notation is dropped for
convenience.
Using dimensionless numbers
\begin{align*}
\Rey &= \frac{\hat\rho \hat{v} \hat{L}}{\hat\mu}
&
\Ca &= \frac{\hat{v} \hat\mu}{\tension}
&
\Eu &= \frac{\hat{p}}{\hat\rho \hat{v}^2}
&
\Fr &= \frac{\hat{v}}{\sqrt{g \hat{L}}}
&
S &= \frac{\hat\mob}{\hat{v}}
~,
\end{align*}
one obtains (see Appendix~\ref{app:non-dim} for details)
\begin{subequations}
\label{eq:non-dim_ns}
\begin{alignat}{2}
\label{eq:non-dim_mass}
\frac{\partial \rho}{\partial t} + \nabla \cdot (\rho \vv) &= 0
& & ~, \text{ in }\Po ~, \\
\label{eq:non-dim_momentum}
\frac{\partial}{\partial t} (\rho \vv)
+ \nabla \cdot (\rho \vv \otimes \vv)
    &=
- \Eu \nabla p
+ \frac{1}{\Rey} \nabla \cdot \left( \mu \left( \nabla \vv + (\nabla \vv)^T
             - \frac{2}{3} (\nabla \cdot \vv)\mathbf{I} \right) \right)
\\ \notag
&\phantom{=}
- \frac{1}{\Fr^2} \rho \mathbf{z}
- \frac{\e}{\Ca} \frac{3\xi}{2}
    \nabla \cdot (\nabla \pf \otimes \nabla \pf)
& & ~, \text{ in }\Po ~,
\end{alignat}
\end{subequations}
with boundary conditions
\begin{equation}
\label{eq:non-dim_bc_slip}
\vv = - \e \slip (\p_\n \vv_\ve{t}) \ve{t} ~,  \text{ on } \InterfaceCell
\end{equation}

For the phase-field equation the non-dimensionalization
yields (see Appendix~\ref{app:non-dim} for details)
\begin{equation}
\label{eq:non-dim_pf}
\frac{\p \pf}{\p t} - S \e^1 \xi \nabla^2 \pf + \nabla \cdot (\vv \pf)
 = -S \e^{-1} \xi^{-1} P'(\pf)
~, \text{ in }\Po ~,
\end{equation}
with boundary condition
\begin{equation}
\label{eq:non-dim_bc_pf}
\p_\n \pf = -\e^{-1} \xi^{-1} \cos(\contactAngle) \sqrt{2P(\pf)} ~,  \text{ on } \InterfaceCell
~.
\end{equation}

\paragraph{Phase-constant density and viscosity}
Using the above equations \eqref{eq:fluid_properties_linear} for density and viscosity and the
reference  values \eqref{eq:non-dim_reference}, the non-dimensional viscosity and density in
\eqref{eq:non-dim_ns} are
\begin{align}
\mu(\pf) =& \frac{1}{\hat\mu}(\mu_1 \pf + \mu_2 (1-\pf))
= 1 + \pf (M - 1)
~,
\\
\rho(\pf) =& 1 + \pf (R - 1)
~.
\end{align}

\subsection{Periodic Homogenization}
\label{sec:homogenization}
Given the scale separation $\e$ we introduce the micro-scale coordinate $\y = \e^{-1} \x$ and assume all
unknowns can be written as an asymptotic expansion in $\e$, depending on both $\x$ and $\y$ with periodicity in the unit cell $Y$.
For an unknown $\psi \in \{p, \vv, \pf\}$ we introduce $Y$-periodic $\psi_k(t, \x, \y)$ such that
\begin{equation*}
\psi(t, \x) = \sum_{k=0}^\infty \e^k \psi_k \left(t, \x, \frac{\x}{\e} \right) ~.
\end{equation*}
The spatial derivatives are rewritten as
\begin{equation*}
\nabla \psi = \nabla_\x \sum_{k=0}^\infty \e^k \psi_k (t, \x, \y) + \frac{1}{\e} \nabla_\y \sum_{k=0}^\infty \e^k \psi_k (t, \x, \y) ~.
\end{equation*}
For the upscaling we consider a flow regime where Darcy's law is considered valid, with laminar flow
driven by the pressure drop and capillary forces, and where advective and diffusive time scales are of the same order. 
For the phase distribution the phase-field diffusivity $\mob$ should be comparable to the
micro-scale advection, captured by $S \in \Oe{0}$.
As will be shown in Remark~\ref{rem:scaling_diff}, with $S\in\Oe{1}$ the advective term separates in
the upscaling process and yields a restrictive pore-scale equation as well as introducing
mixed-scale derivatives into the phase-field equation.
Choosing non-dimensional numbers $\Ca \in \Oe{0}$, $\Rey \in \Oe{0}$, $\Eu \in \Oe{-2}$ and
$\Fr \in \Oe{0}$ the leading terms for $\e \to 0$ are given as follows (see
Appendix~\ref{app:homogenization} for details).

For the phase dependent fluid properties we observe
\begin{subequations}
\begin{align}
\label{eq:hom_viscosity}
\mu(\pf) &= (1+\pf(M-1))
    = \underbrace{1 + \pf_0(M-1)}_{=\mu(\pf_0)} + \e \pf_1(M-1) + \Oe{2}
\\
\label{eq:hom_density}
\rho(\pf) &= (1+\pf(R-1))
    = \underbrace{1 + \pf_0(R-1)}_{=\rho(\pf_0)}
    + \e \underbrace{\pf_1(R-1)}_{=\rho_1} + \Oe{2}
~,
\intertext{and for the double-well potential of the phase-field equation, due to its polynomial
    structure,}
P'(\pf) &= P'(\pf_0) + \e \pf_1 P''(\pf_0) + \Oe{2}
~.
\end{align}
\end{subequations}

Denoting $\overline{\Eu} := \e^2 \Eu$, we have for \eqref{eq:non-dim_mass},
\eqref{eq:non-dim_momentum} and \eqref{eq:non-dim_pf} in $\Po$
\begin{subequations}
\begin{align}
\Oe{1} =& \e^{-1} \nabla_\y \cdot (\rho(\pf_0) \vv_0)
\label{eq:expansions_mass_-1}
\\
&+ \e^{0} \bigg[ \frac{\partial\rho(\pf_0)}{\partial t}
        + \nabla_\y \cdot (\rho(\pf_0) \vv_1 + \rho_1 \vv_0)
        + \nabla_\x \cdot (\rho(\pf_0) \vv_0) \bigg]
\notag
\\
\Oe{-1} =& \e^{-3} \left[ -\overline\Eu \nabla_\y p_0 \right]
\label{eq:expansions_mom_-3}
\\
&+ \e^{-2} \bigg[ -\overline\Eu (\nabla_x p_0 + \nabla_\y p_1)
\notag\\&
+ \frac{1}{\Rey} \nabla_\y \cdot \left( \mu(\pf_0) \left( \nabla_\y \vv_0 +(\nabla_\y \vv_0)^T - \frac{2}{3}
            (\nabla_\y \cdot \vv_0) I \right)\right)
\notag
\\
& - \frac{1}{\Ca} \frac{3\xi}{2} \nabla_\y \cdot (\nabla_\y \pf_{0} \otimes \nabla_\y \pf_{0} ) \bigg]
\notag
\\
\Oe{1} =& \e^{-1} [ \nabla_\y \cdot (\vv_0 \pf_{0}) - S \xi \nabla_\y^2 \pf_{0}
+ S \xi^{-1} P'(\pf_{0})]
\label{eq:expansions_pf_-1}
\\
&+ \e^{0} \bigg[ \frac{\p \pf_{0}}{\p t} + \nabla_\y \cdot
(\vv_1 \pf_{0} + \vv_0 \pf_{1}) + \nabla_\x \cdot (\vv_0 \pf_{0})
+S \xi^{-1} P''(\pf_{0}) \pf_{1}
\notag\\
&-S \xi ( \nabla_\y \cdot \nabla_\x \pf_{0} + \nabla_\x \cdot \nabla_\y \pf_{0}
 + \nabla_\y^2 \pf_{0} )
\bigg]
\notag
\intertext{and for boundary conditions
\eqref{eq:non-dim_bc_slip} and
\eqref{eq:non-dim_bc_pf} on $\InterfaceCell$}
\label{eq:expansions_bc_slip}
\Oe{1} &= \e^{0} \Big[ \vv_0 + \slip (\nabla_\y \vv_{0,\ve{t}} \cdot \n) \ve{t} \Big]
\\
\label{eq:expansions_bc_pf}
\Oe{0} &= \e^{-1} \Big[ \nabla_\y \pf_0 \cdot \n  + \xi^{-1} \cos(\contactAngle) \sqrt{2P(\pf_0)} \Big]
~.
\end{align}
\end{subequations}

\subsubsection{Phase field}
From the leading order term of \eqref{eq:expansions_pf_-1} we therefore obtain a stationary
equation involving only local derivatives.
\begin{equation}
\label{eq:cell_pf}
- S \xi \nabla_\y^2 \pf_{0}
+ \nabla_\y \cdot (\vv_0 \pf_{0})
= S \xi^{-1} P'(\pf_{0})
~, \qquad \text{ in }\Po
~,
\end{equation}
together with boundary condition \eqref{eq:expansions_bc_pf}
\begin{equation}
\label{eq:cell_pf_bc}
\nabla_\y \pf_0 \cdot \n = \pfBCcell ~,
\qquad\text{ on }\InterfaceCell~.
\end{equation}
\begin{remark}
\label{rem:scaling_diff}
If the phase-field equation is dominated by advection ($S \in \Oe{1}$)
the leading order advective term of \eqref{eq:expansions_bc_pf} would be isolated at $\Oe{-2}$ as
\begin{equation}
\label{eq:dominant_advection}
\nabla_\y \cdot ( \vv_0 \pf_{0} ) = 0
~.
\end{equation}
Applied to the leading order term in \eqref{eq:expansions_mass_-1}
this yields a divergence free velocity field and reduces to
\begin{equation}
\label{eq:restriction_flowlines}
\vv_0 \cdot \nabla_\y \pf_{0} = 0
~,
\end{equation}
resulting in a strong limitation on modeled problems.
At the order $\Oe{-1}$ the equation \eqref{eq:expansions_pf_-1} would instead contain the advective
terms of the next order, yielding
\begin{equation}
\label{eq:mixed_derivatives}
0 = 
\nabla_\y \cdot (\vv_1 \pf_0 + \vv_0 \pf_1) + \nabla_\x \cdot (\vv_0 \pf_0)
- S \xi \nabla_\y^2 \pf_{0} + S \xi^{-1} P'(\pf_{0})
\end{equation}
with an undesireable mix of scales that leads to neither a pore-scale cell problem nor an
effective equation for the Darcy scale.
A balance of the two terms is expressed by $S \in \mathcal{O}(1)$, which yields
\eqref{eq:cell_pf} instead of \eqref{eq:dominant_advection} as the leading order
equation.
\end{remark}
This local phase-field equation \eqref{eq:cell_pf} is underconstrained, admitting among others the trivial solutions of
$\pf_0=0$ and $\pf_1=1$ for divergence free velocity fields.
The local cell problem does not offer a way to compute the saturation of fluid 1 as the mean
integral of $\pf_0$ but rather requires it as a constraint as in \cite{sharmin2022}.

From the next order $\Oe{0}$ terms of the phase-field equation \eqref{eq:expansions_pf_-1} we
obtain
\begin{equation}
\label{eq:expansion_pf_0}
\begin{aligned}
0&=
\frac{\partial \pf_0}{\partial t}
- S\xi( \nabla_\y^2 \pf_{1}
        + \nabla_\y \cdot (\nabla_\x \pf_{0})
        + \nabla_\x \cdot (\nabla_\y \pf_{0})
        )
\\&\phantom{=}
+ ( \nabla_\y \cdot (\vv_1\pf_{0} + \vv_0\pf_{1})
        + \nabla_\x \cdot (\vv_0 \pf_{0})
        )
+ S \xi^{-1} P''(\pf_{0}) \pf_{1}
~.
\end{aligned}
\end{equation}
We would like to use this instationary equation to update the saturation to in turn constrain
the stationary phase-field equation obtained from the leading order term of the asymptotic expansion.
After integrating over the constant local periodicity cell $\Po$ and using the periodicity with
respect to $\y$, \eqref{eq:expansion_pf_0} reduces to
\begin{equation}
\label{eq:pf_0_int}
\frac{\partial}{\partial t} \int_\Po \pf_{0} \dd\y
+ \nabla_\x \cdot \int_\Po \vv_0 \pf_{0} \dd\y
+ \int_\Po \frac{S}{\xi} P''(\pf_{0}) \pf_{1} \dd\y
=0
~.
\end{equation}
While the first two terms yield a macroscopic equation in saturation and phase-specific velocity,
the last term includes the additional unknown $\pf_{1}$.
Assuming a homogeneous porous medium, with $\Po$ not depending on $\x$, one can use a
solvability constraint to show that this term is zero and obtain the desired saturation equation.

Integrating the leading order terms \eqref{eq:cell_pf} and applying the divergence theorem and periodicity
conditions, one is left with only the potential derivative.
\begin{equation}
\label{eq:p'_int}
0 = \underbrace{\int_\Po -S \xi \nabla_\y \cdot (\nabla_\y \pf_0) + \nabla_\y \cdot (\vv_0 \pf_0) \dd\y
}_{=0}
+ \int_\Po \frac{S}{\xi} P'(\pf_0) \dd\y
\end{equation}
One now views the third term in \eqref{eq:pf_0_int}
as a Fredholm operator of index zero $\mathcal{L}(\pf_0)$ applied to $\pf_1$. Applying it
instead to $\nabla_\x \pf_0$, using the chain rule
and the assumption of $\Po$ not depending on $\x$, we obtain
\begin{equation}
\mathcal{L}(\pf_0) (\nabla_\x \pf_0) = \int_\Po \frac{S}{\xi} P''(\pf_0) \nabla_\x \pf_0 \dd\y
= \int_\Po \frac{S}{\xi} \nabla_\x (P'(\pf_0)) \dd\y
= \frac{S}{\xi} \nabla_\x \int_\Po P'(\pf_0) \dd\y
~.
\end{equation}
Together with the previously derived information about $P'$ from \eqref{eq:p'_int}, one sees that
$\nabla_\x \pf_0$ is an
element of the kernel of $\mathcal{L}(\pf_0)$.
Rewriting \eqref{eq:pf_0_int} as $-\mathcal{L}(\pf_0)\pf_1 = A(\pf_0, \vv_0)$ with
\begin{equation}
A(\pf_0, \vv_0) =
\int_\Po \frac{\partial \pf_0}{\partial t} + \nabla_\x \cdot (\vv_0 \pf_0) \dd\y ~,
\end{equation}
this yields the solvability constraint
\begin{equation}
0 = <A(\pf_0, \vv_0), \nabla_\x \pf_0>
= \int_\Omega \nabla_\x \pf_0 (\x,\y) \left(\int_\Po \frac{\partial \pf_0}{\partial t}
 + \nabla_\x \cdot (\vv_0 \pf_0) \dd\y
 \right) (\x) \dd\x
~.
\end{equation}
To avoid trivial behavior ($\nabla_\x \pf_0 \equiv 0$, $\nabla_\x \pf_0$ independent of
$\x$), the integral over the local pore-space $\Po$ must disappear.
Using again the assumption of a
stationary and homogeneous pore space $\Po$ with constant porosity $\porosity$,
this yields the saturation equation
\begin{equation}
\label{eq:saturation_constraint_0}
\frac{\partial}{\partial t} \int_\Po \pf_0 \dd\y + \nabla_\x \cdot \int_\Po \vv_0 \pf_0 \dd\y = 0
~.
\end{equation}
Introducing averaged quantities for the saturation of fluid 1, 
$S^{(1)} = \poro^{-1} \int_\Po \pf_0 \dd\y$, as well as velocities
\begin{equation}
\bar\vv^{(1)} = \poro^{-1}\int_\Po \pf_0 \vv_0 \dd\y
~,\qquad
\bar\vv^{(2)} = \poro^{-1}\int_\Po (1 - \pf_0) \vv_0 \dd\y
~,\qquad
\bar\vv = \poro^{-1}\int_\Po \vv_0 \dd\y
~,
\end{equation}
we obtain
\begin{equation}
\label{eq:saturation_1}
\frac{\partial}{\partial t}
    \int_\Po \pf_0 \dd\y
+ \nabla_\x \cdot
    \int_{\Po} \vv_0 \pf_0\dd\y
= \porosity\left(
        \frac{\partial}{\partial t} S^{(1)}(\x) + \nabla_\x \cdot \bar\vv^{(1)}(\x)
    \right)
= 0
~.
\end{equation}

The second order term of the mass conservation equation \eqref{eq:expansions_mass_-1}
\begin{equation*}
0= \frac{\partial\rho(\pf_0)}{\partial t}
+ \nabla_\y \cdot (\rho(\pf_0) \vv_1 + \rho_1 \vv_0)
+ \nabla_\x \cdot (\rho(\pf_0) \vv_0)
~,
\end{equation*}
yields, using \eqref{eq:hom_density} and after integration over $\Po$,
\begin{equation}
\label{eq:mass_averaged:saturation}
(R-1) \frac{\partial}{\partial t} \int_\Po \pf_{0} \dd\y + (R-1) \nabla_\x \cdot \int_\Po \pf_{0} \vv_0 \dd\y
+ \nabla_\x \cdot \int_\Po \vv_0 \dd\y ~,
\end{equation}
a second conservation equation.

Inserting \eqref{eq:saturation_1} into \eqref{eq:mass_averaged:saturation} yields
$\nabla_\x \cdot \bar{\vv} = 0$ and with $S^{(2)} = 1 - S^{(1)}$  a saturation equation for the second phase, corresponding to $\pf_0 = 0$,
\begin{equation}
\label{eq:saturation_2}
0 = \frac{\partial}{\partial t} (S^{(1)} + S^{(2)})
= \frac{\partial}{\partial t} S^{(2)}
+ \frac{1}{\porosity} \nabla_\x \cdot \int_\Po (1-\pf_0) \vv_0 \dd\y
= \frac{\partial}{\partial t} S^{(2)} + \nabla_\x \cdot \bar\vv^{(2)}
~.
\end{equation}

These macroscopic equations \eqref{eq:saturation_1}, \eqref{eq:saturation_2}, through the saturation,
yield an integral constraint for the local
phase-field. Together with the stationary equation \eqref{eq:cell_pf} and boundary condition
\eqref{eq:cell_pf_bc} we
obtain the cell problem
\begin{equation}
  \label{eq:cell_phasefield}
  \begin{cases}
    \nabla_\y \cdot (\vv_0 \pf_0) = \mb \xi \nabla_\y^2 \pf_0 - \mb \xi^{-1} P'(\pf_0)
    & \text{ in } \Po
    ~,\\
    \nabla_\y \pf_0 \cdot \n = \pfBCcell
    & \text{ on } \InterfaceCell
    ~,\\
    \pf_0 \text{ is $Y$-periodic and } \int_\Po \pf_0 = \poro S^{(1)} ~.
  \end{cases}
\end{equation}
While this ensures mass-conservation it does not fully
prescribe a distribution of the phases. The stationary phase-field equation
\eqref{eq:cell_phasefield} still defines the
profile of the diffuse transition zone but the position is not clear. This phase distribution
however is central to the equations for fluid flow.
To obtain a meaningful solution to this stationary problem one could introduce an artifical time
evolution, taking care to introduce additional source terms to the non-conservative Allen-Cahn
equation in order to enforce the saturation constraint.
That is, one would instead solve
\begin{equation}
  \label{eq:cell_phasefield_instationary}
  \begin{cases}
    \frac{\partial \pf_0}{\partial \tau} +
    \nabla_\y \cdot (\vv_0 \pf_0) = \mb \xi \nabla_\y^2 \pf_0 - \mb \xi^{-1} P'(\pf_0)
    &\\\phantom{\frac{\partial \pf_0}{\partial \tau} + \nabla_\y \cdot (\vv_0 \pf_0) =}
    + \frac{\mb\xi^{-1}}{\porosity} \int_\Po P'(\pf_0) \dd\y
    - \interface \left( \int_\Po \pf_0 \dd\y - S^{(1)} \right)
    & \text{ in } \Po
    ~,\\
    \nabla_\y \pf_0 \cdot \n = \pfBCcell
    & \text{ on } \InterfaceCell
    ~,\\
    \pf_0 \text{ is $Y$-periodic.}
  \end{cases}
\end{equation}
The third term on the right-hand side is introduced in order to make the Allen-Cahn equation conservative 
\cite{bringedal2020con}. The
fourth term on the right-hand side moves the system towards the desired saturation and can be localized to the
interface using an interface indicator function
\begin{equation*}
\interface(\pf_0) = 4 \xi^{-1} \pf_0 (1-\pf_0) ~.
\end{equation*}

\subsubsection{Fluid flow}\label{sec:momentum_dilatation}
The derivation of the local cell problems for the fluid velocity follows the common approach of
separating derivatives of the two scales and using the linear structures of the equations to obtain
different problems for the effect of macroscopic pressure gradients, and additionally one to capture the contribution of the surface tension at the interface.

From the leading order term of the momentum equation \eqref{eq:expansions_mom_-3} we obtain
$\nabla_\y p_0 = 0$ and from the next order term we get
\begin{equation*}
\begin{aligned}
0 =&
 -\overline\Eu (\nabla_\x p_0 + \nabla_\y p_1)
+ \frac{1}{\Rey} \nabla_\y \cdot \left( \mu(\pf_0) \left( \nabla_\y \vv_0 +(\nabla_\y \vv_0)^T
            - \frac{2}{3} (\nabla_\y \cdot \vv_0)\tens{I} \right)\right)
\\
 &- \frac{1}{\Ca} \frac{3\xi}{2} \nabla_\y \cdot (\nabla_\y \pf_{0} \otimes \nabla_\y \pf_{0} )
\notag
\end{aligned}
\end{equation*}
or equivalently, using $\mu(\pf_0(t, \x, \y)) = 1 + \pf_0(t, \x, \y) (M-1)$,
\begin{equation}
\label{eq:cell_mom_dilatation}
\begin{aligned}
&\overline\Eu \nabla_\y p_1
- \frac{M-1}{\Rey} \nabla_\y \pf_{0} \cdot \left( \nabla_\y \vv_0 +(\nabla_\y \vv_0)^T
        -\frac{2}{3} (\nabla_\y \cdot \vv_0) \tens{I} \right)
\\
&- \frac{\mu(\pf_0)}{\Rey} \left( \nabla_\y^2 \vv_0 - \frac{1}{3} \nabla_\y(\nabla_\y \cdot \vv_0) \right)
\\
=&
 -\overline\Eu \nabla_\x p_0
 - \frac{1}{\Ca} \frac{3\xi}{2} \nabla_\y \cdot (\nabla_\y \pf_{0} \otimes \nabla_\y \pf_{0} )
~.
\end{aligned}
\end{equation}
Due to the linear structure in the unknowns $p_1$ and $\vv_0$ as well as the contributions of
$p_0(t, \x)$ and $\pf_{0}(t, \x, \y)$, we can find $Y$-periodic functions $\vw_j(t, \x, \y)$,
$\Pi_j(t, \x, \y)$,
$j=0,\ldots, \dimension$ such that
\begin{align}
\label{eq:velocity_decomposition}
\vv_0 =& -\sum_{j=1}^\dimension \vw_j \partial_{\x_j} p_0 - \vw_0
~, \\
p_1 =& \sum_{j=1}^\dimension \Pi_j \partial_{\x_j} p_0 + \Pi_0
+ \tilde{p}_1(t, \x)
~,
\end{align}
with a function $\tilde{p}_1(t, \x)$ independent of $\y$ and thus not relevant for
\eqref{eq:cell_mom_dilatation}.
These velocity contributions and local pressures can be obtained from the following auxiliary cell
problems, where we denote the symmetrized
gradient $2\varepsilon_\y(\vw) = \nabla_\y \vw + (\nabla_\y \vw)^T$,
\newcommand{\flowBCcell}[1]{-\slip (\partial_\n \vw_{#1\ve{t}}) \ve{t}}
\begin{equation}\label{eq:velocity_pressure}
\begin{cases}
\overline\Eu (e_j + \nabla_\y \Pi_j)
=
\frac{1}{\Rey} \nabla_\y \cdot ( \mu(\pf_{0}) 2\varepsilon_\y(\vw_j)
        - \frac{2}{3} (\nabla_\y \cdot \vw_j)\tens{I})
, & \text{ in } \mathcal{P} ~,
\\
\nabla_\y \cdot (\rho(\pf_{0}) \vw_j) = 0
, & \text{ in } \mathcal{P} ~,
\\
\vw_j = \flowBCcell{j,}
, & \text{ on } \InterfaceCell ~,
\\
\Pi_j, \vw_j \text{ are } Y\text{-periodic and} \int_\mathcal{P} \Pi_j \dd\y = 0
,
\end{cases}
\end{equation}
for $j \in \{ 1, \ldots , \dimension \}$ and
\begin{equation}\label{eq:velocity_surface}
\begin{cases}
\overline\Eu \nabla_\y \Pi_0
=
\frac{1}{\Rey} \nabla_\y \cdot ( \mu(\pf_{0}) 2\varepsilon_\y(\vw_0)
        -\frac{2}{3} (\nabla_\y \cdot \vw_0) \tens{I} )
\\
\hphantom{\overline\Eu \nabla_\y \Pi_0 =}
- \frac{1}{\Ca} \frac{3\xi}{2} \nabla_\y \cdot (\nabla_\y \pf_{0} \otimes \nabla_\y \pf_{0} )
, & \text{ in } \mathcal{P} ~,
\\
\nabla_\y \cdot (\rho(\pf_{0})\vw_0) = 0
, & \text{ in } \mathcal{P} ~,
\\
\vw_0 = \flowBCcell{0,}
, & \text{ on } \InterfaceCell ~,
\\
\Pi_0, \vw_0 \text{ are } Y\text{-periodic and} \int_\mathcal{P} \Pi_0 \dd\y = 0
,
\end{cases}
\end{equation}

From these solutions we define tensors capturing the effective behavior of the fluid, denoting the
phase indicator of phase $k$ as
\begin{equation}
\pf^{(k)} = \begin{cases} \pf_0  ,\,&k=1 \\ 1 - \pf_0  ,\,&k=2 \end{cases}
\end{equation}
and the components of $\vw_j$ as $\vw_{j,i}$.
\begin{equation}
\effPerm^{(k)}_{ij} := \frac{1}{\porosity}\int_\Po \pf^{(k)} \vw_{j,i} \dd\y
\qquad
\effSurf^{(k)}_{i} := \frac{1}{\porosity}\int_\Po \pf^{(k)} \vw_{0,i} \dd\y
\end{equation}
Multiplying \eqref{eq:velocity_decomposition} with $\pf^{(k)}$ and integrating over $\Po$, we obtain the
macroscopic velocity equations containing the effective parameters.
\begin{equation}
\begin{aligned}
\bar{\vv}^{(k)} = \varphi^{-1} \int_\Po \pf^{(k)} \vv_0 \dd\y
&= \varphi^{-1} \int_\Po \pf^{(k)} (-\sum_{j=1}^\dimension \vw_j \partial_{\x_j} p_0 - \vw_0) \dd\y
\\
&= - \sum_{j=1}^\dimension \left(\varphi^{-1} \int_\Po \pf^{(k)} \vw_j \dd\y\right) (\partial_{\x_j} p_0)
    - \varphi^{-1} \int_\Po \pf^{(k)} \vw_0 \dd\y
\\
&= - \effPerm^{(k)} \nabla_x p_0 - \effSurf^{(k)}
~.
\end{aligned}
\end{equation}

\subsection{Two-scale model}
In summary we arrive at a micro-macro model, coupling a Darcy-scale flow problem reminiscent of the
extended two-phase Darcy's law to $d+2$ pore-scale cell-problems on $Y = [0,1]^d$ at every
point of the domain.
We solve for five macroscopic unknowns $p_0(t,\x)$, $S^{(1)}(t,\x)$, $S^{(2)}(t,\x)$,
$\bar\vv^{(1)}(t, \x)$, $\bar\vv^{(2)}(t, \x)$ with
$1 = S^{(1)} + S^{(2)}$, and for every set of cell problems the local unknowns $\pf_0$ and $\vw_i$,
$i = 0, \ldots \dimension$.

The macro-scale model reminds of Darcy's law with one shared pressure unknown as well as
additional effective parameters $\effSurf^{(k)}$ which, just as $\effPerm^{(k)}$, are computed
from cell-problem variables.
For convenience we drop the subscript $\x$ on the macroscopic gradient.
\begin{subequations}
\begin{align}
  \label{eq:twoscale_sat}
  0 &= \frac{\partial}{\partial t} S^{(k)} + \nabla \cdot \bar\vv^{(k)} \quad  \text{, in } \Omega \text{, } k=1,2
  \\
  \bar\vv^{(k)} &= -\effPerm^{(k)} \nabla p - \effSurf^{(k)} \quad \text{, in } \Omega \text{, } k=1,2
  \label{eq:twoscale_vel}\\
  \effPerm^{(k)}_{ij} &= \frac{1}{|\Po|} \int_\Po \pf^{(k)} (\vw_j)_i \dd\y
  ~,\qquad 
  \label{eq:twoscale_eff}
  \effSurf^{(k)}_i = \frac{1}{|\Po|} \int_\Po \pf^{(k)} (\vw_0)_i \dd\y
\end{align}
\end{subequations}
We see that the effective parameters $\effPerm^{(k)}$ represent the phase-specific effective mobilities,
containing information about both the absolute permeability of the porous medium and the
interactions between the two fluids.
Both contributions are in general anisotropic and can't easily be separated.
For isotropic geometries the intrinsic permeability $\kappa_\mathrm{abs}$ is a scalar and the effective mobility can be
written as $\kappa_\mathrm{abs} \effPerm^{(k)}_\mathrm{rel}(\pf_0) /\mu_k$ with the relative permeability
$\effPerm^{(k)}_\mathrm{rel} \in \mathbb{R}^{\dimension \times \dimension}$ 
depending on the local phase distribution $\pf_0$ through the cell problems
\eqref{eq:two_scale_cell_permeability}.
The second effective parameter $\effSurf^{(k)}$ captures effective flow due to surface tension forces
between the two fluids.

In return, the cell problems depend on the local value of the macroscopic saturation through a
constraint on the phase-field function.
\begin{equation}
  \label{eq:two_scale_cell_phasefield}
  \begin{cases}
    \nabla_\y \cdot (\vv_0 \pf_0) = \mb \xi \nabla_\y^2 \pf_0 - \mb \xi^{-1} P'(\pf_0)
    & \text{ in } \Po
    ~,\\
    \nabla_\y \pf_0 \cdot \n = \pfBCcell
    & \text{ on } \InterfaceCell
    ~,\\
    \pf_0 \text{ is $Y$-periodic and } \int_\Po \pf_0 \dd\y = S^{(1)} ~,
  \end{cases}
\end{equation}
with $\vv_0$ depending on the global pressure gradient and local flow functions $\vw_j$ as defined
in \eqref{eq:velocity_decomposition}.
These in turn depend on the phase field and are computed through the cell problems defined in
\eqref{eq:velocity_pressure} and \eqref{eq:velocity_surface}
\begin{equation}
\label{eq:two_scale_cell_permeability}
\begin{cases}
\overline\Eu (e_j + \nabla_\y \Pi_j)
=
\frac{1}{\Rey} \nabla_\y \cdot ( \mu(\pf_0) 2\varepsilon_\y(\vw_j)
        - \frac{2}{3} (\nabla_\y \cdot \vw_j)\tens{I})
, & \text{ in } \mathcal{P} ~,
\\
\nabla_\y \cdot (\rho(\pf_0) \vw_j) = 0
, & \text{ in } \mathcal{P} ~,
\\
\vw_j = \flowBCcell{j,}
, & \text{ on } \InterfaceCell ~,
\\
\Pi_j, \vw_j \text{ are } Y\text{-periodic and} \int_\mathcal{P} \Pi_j \dd\y = 0
,
\end{cases}
\end{equation}
for $j \in \{ 1, \ldots , d \}$ and
\begin{equation}
\label{eq:two_scale_cell_interface}
\begin{cases}
\overline\Eu \nabla_\y \Pi_0
=
\frac{1}{\Rey} \nabla_\y \cdot ( \mu(\pf_0) 2\varepsilon_\y(\vw_0)
        -\frac{2}{3} (\nabla_\y \cdot \vw_0) \tens{I} )
\\
\hphantom{\overline\Eu \nabla_\y \Pi_0 =}
- \frac{1}{\Ca} \frac{3\xi}{2} \nabla_\y \cdot (\nabla_\y \pf_0 \otimes \nabla_\y \pf_0 )
, & \text{ in } \mathcal{P} ~,
\\
\nabla_\y \cdot (\rho(\pf_0)\vw_0) = 0
, & \text{ in } \mathcal{P} ~,
\\
\vw_0 = \flowBCcell{0,}
, & \text{ on } \InterfaceCell ~,
\\
\Pi_0, \vw_0 \text{ are } Y\text{-periodic and} \int_\mathcal{P} \Pi_0 \dd\y = 0
.
\end{cases}
\end{equation}
\section{Numeric Investigation}\label{sec:numerics}
In the following we present a simulation using the tightly coupled pore-scale model from
Section~\ref{sec:pf_model}, as well as the investigation of
the cell problems for velocity components $\vw_j$ and the behavior of the computed
effective parameters.
For the former we simulate the movement of two fluids through a channel geometry modeling a
single pore.
For the latter we consider two exemplary geometries, denoted "cross" and "obstacle", and investigate
the relative permeability saturation relation for different fluid properties and fluid
distributions.

\subsection{Implementation}
\begin{figure}
\center
\includegraphics[width=0.5\textwidth]{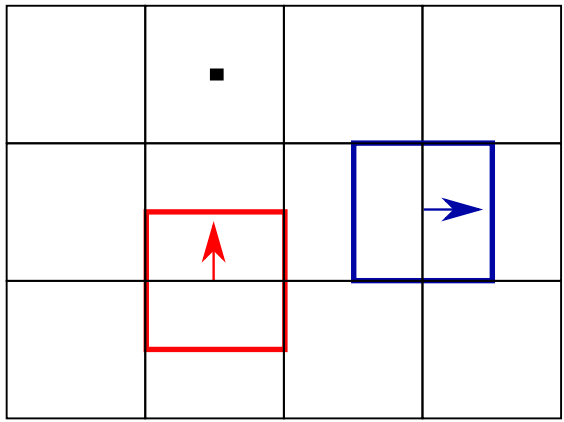}
\caption{Staggered discretization}
\label{fig:staggered}
\end{figure}
For the numerical investigation
all model equations are implemented in \dumux \cite{KOCH2021423}
, using Newton's method as a nonlinear solver and backward Euler for temporal discretization in
transient problems.
The equations are discretized in space with the finite volume method on a regular rectangular
grid, using Dune-SPGrid \cite{nolte2011} and an adapted Dune-Subgrid \cite{graeser2009} to model the considered geometries with periodic
boundaries.
For the phase-field and pressure unknowns a cell-centered approach is used, with control volumes
equal to grid cells and degrees of freedom placed at their centers.
Fluxes over control-volume faces are approximated using the two adjacent values.
For the velocities a staggered discretization is used (see Figure~\ref{fig:staggered}), with
separate degrees of freedom for each velocity component placed at the edges
between grid cells and control volumes centered around them.

While the cell problem for the phase field uses only the cell-centered discretization,
the pore-scale problem and 
the cell problems for velocity contributions require a coupled system of
cell-centered pressures and staggered velocities.
In \dumux this is achieved using a multidomain formulation with a coupling manager handling the
volume coupling of the different discretizations of a shared domain \cite{KOCH2021423}.
For the pore-scale model the existing coupling manager for Navier-Stokes models was extended to
additionally exchange phase-field data.

For the phase-field equation a boundary condition prescribing the flux implements the contact angle
condition.
In the velocity problems a no-flux condition for the mass conservation and homogeneous Dirichlet
conditions for the normal velocity are used at the fluid-solid interfaces.
To implement Navier slip boundary conditions a solution-dependent flux can be described.
In the case of the cell problem capturing surface tension effects the additional flux can be given
based on the prescribed contact angle.
On the outer boundary of the reference cell $Y$ periodicity conditions are prescribed.

\subsection{Pore-scale simulation}
We present an examplary simulation of the fully coupled and non-dimensionalized pore-scale model from
Section~\ref{sec:non-dim}.
In a channel geometry containing both phases the system conforms to a prescribed contact angle and
the interface is advected by a pressure gradient.
This could be applied to simulating two-phase flow through a narrowing in a porous medium.

\begin{figure}
  \includegraphics[width=\textwidth]{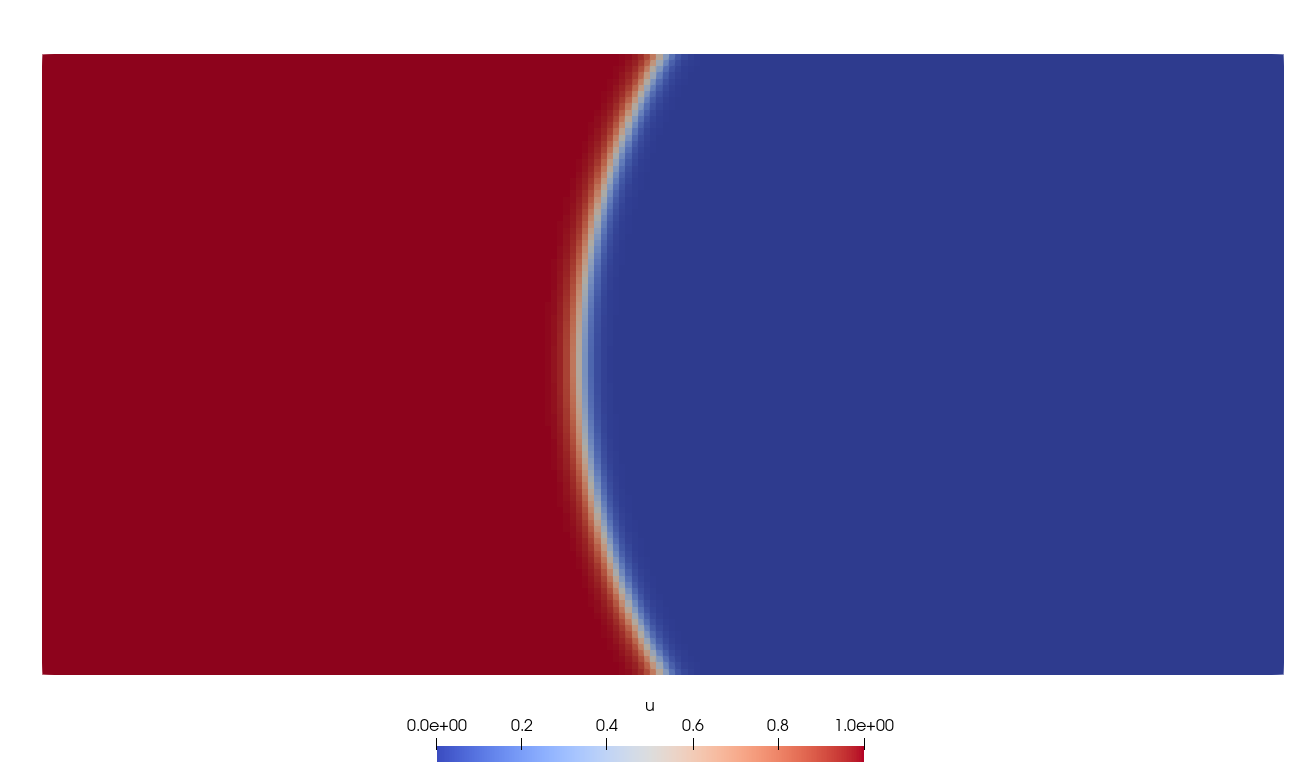}
  \caption{Initial condition of pore-scale simulation of two-phase flow through a channel. $\pf = 1$ (red) corresponds to
  fluid 1, $\pf = 0$ (blue) to fluid 2.}
  \label{fig:pore_scale_sim_initial}
\end{figure}

We consider a two-dimensional domain $\Omega = (0,0.2)\times(0,0.1)$.
At the inlet ($x=0$) and outlet
($x=0.2$) we prescribe fixed pressures and fluxes for the phase-field equation.
The top and bottom boundaries are walls with slip boundary conditions for the velocity and no-flux
boundary conditions for the density. For the phase field we apply the contact angle boundary
condition \eqref{eq:non-dim_bc_pf} with contact angle $\contactAngle = \pi/3$.

We initialize the fluid distribution with fluid 1 on the left and fluid 2 on the right with a curved
and diffuse interface between them, see Figure~\ref{fig:pore_scale_sim_initial}. We prescribe a
non-dimensional slip length of 0.01 and a surface tension of $10^{-6} \mathrm{m^2}$.
The simulation was run for constant fluid properties, using a density of $10^2 \mathrm{kg/m^3}$ and dynamic viscosity of
$10^{-2} \mathrm{kg/(m \cdot s)}$.
Applying a pressure drop of $450 \mathrm{kg/(m \cdot s^2)}$ leads to a maximum flow velocity of
$5.7 \cdot 10^{-4} \mathrm{m/s}$.

\begin{figure}
  \includegraphics[width=\textwidth]{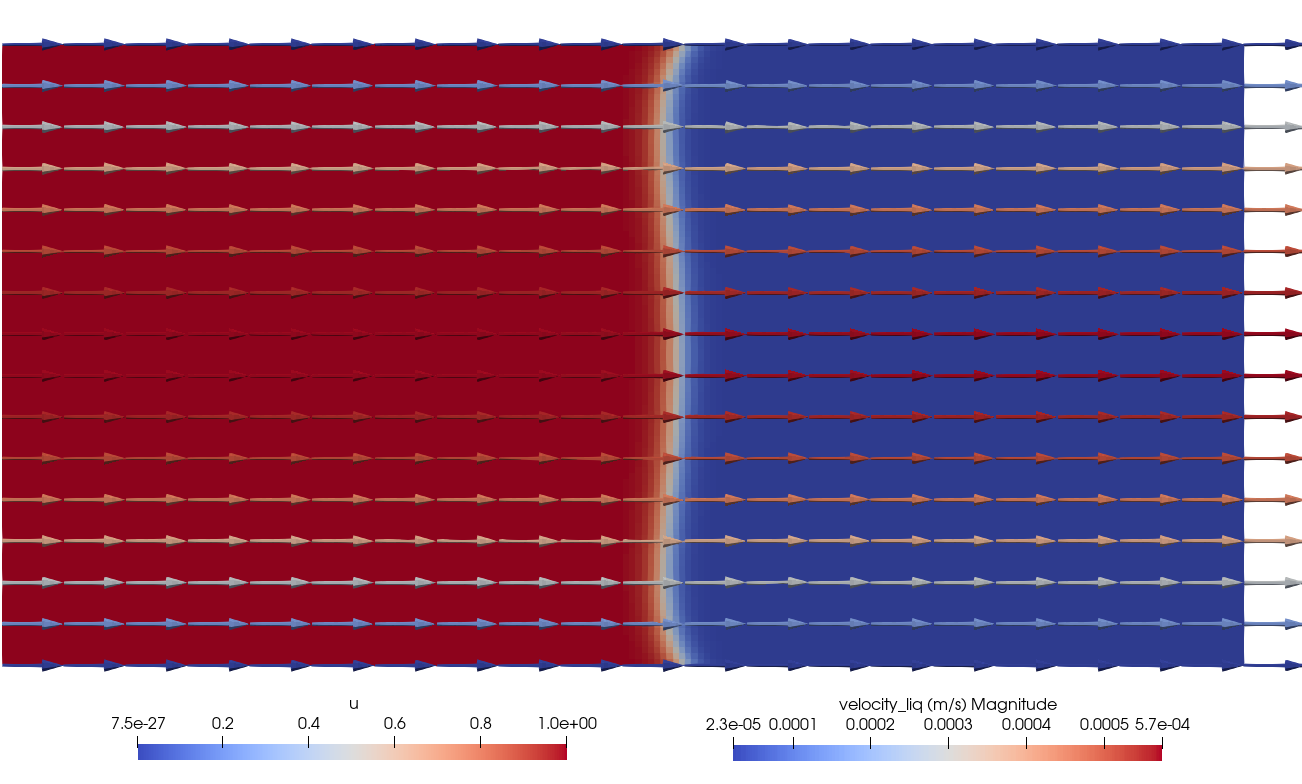}
  \caption{Phase distribution at $t=40\mathrm{s}$ for the pore-scale simulation, overlaid with the
      velocity field. $\pf = 1$ (red) corresponds to
  fluid 1, $\pf = 0$ (blue) to fluid 2.}
  \label{fig:pore_scale_sim}
\end{figure}

Figure~\ref{fig:pore_scale_sim} shows the fluid distribution at $t=40\mathrm{s}$.
Due to very weak surface tension forces in this setup, the velocity field does not deviate
significantly from a parabolic flow profile.
In the center of the channel the interface is advected by these higher velocities, leading to an
inversion of the curvature.
At the three-phase contact point the prescribed contact angle is successfully maintained despite the
near-parabolic profile.

\subsection{Cell problems}
\label{sec:numerics_cell}
We consider the cell problems \eqref{eq:two_scale_cell_permeability} and
\eqref{eq:two_scale_cell_interface} for the velocity components $\vw_j$ and investigate the behavior
of the resulting effective parameters.
The derived two-scale model is similar to two-phase Darcy's law, with effective parameters computed
from cell problems rather than prescribed relative permeability saturation curves.
We investigate how the computed parameters compare to commonly used curves and under which
conditions core assumptions such as a monotone relation is violated.
Anisotropic permeabilities are not the focus of this investigation and we choose simple geometries
such that the intrinsic permeability of the pore geometry can be separated from fluid-fluid interactions.
This allows us to isolate the influence of the fluid distribution on the relative permeability in a
simple manner.

For two different geometries we vary the local phase field
to determine the effects of saturation and phase distribution and the impact of the dimensionless
ratios for viscosity and density as well as the surface tension.
Both the slip length and the contact angle influence the results, see \eg the investigation of dynamic
contact angles in \cite{lunowa2021}, but will not be considered here.
Instead these effects will be controlled by a no-slip condition for the velocity
components and a homogeneous neumann condition for the phase field, encoding a neutral contact
angle.
Since our goal is to address behavior of the effective parameters, we do not couple the flow to the cell
problem for the phase distribution \eqref{eq:two_scale_cell_phasefield}.
Instead we manually prescribe a phase-field function corresponding to a certain saturation and run a few
timesteps of the instationary version of the phase-field equation
\eqref{eq:two_scale_cell_phasefield} without advection in order to enforce the neutral contact angle.
This fluid distribution is then used to solve \eqref{eq:two_scale_cell_permeability} and
\eqref{eq:two_scale_cell_interface}.

\begin{figure}
\begin{subfigure}{0.5\textwidth}
\includegraphics[width=\textwidth]{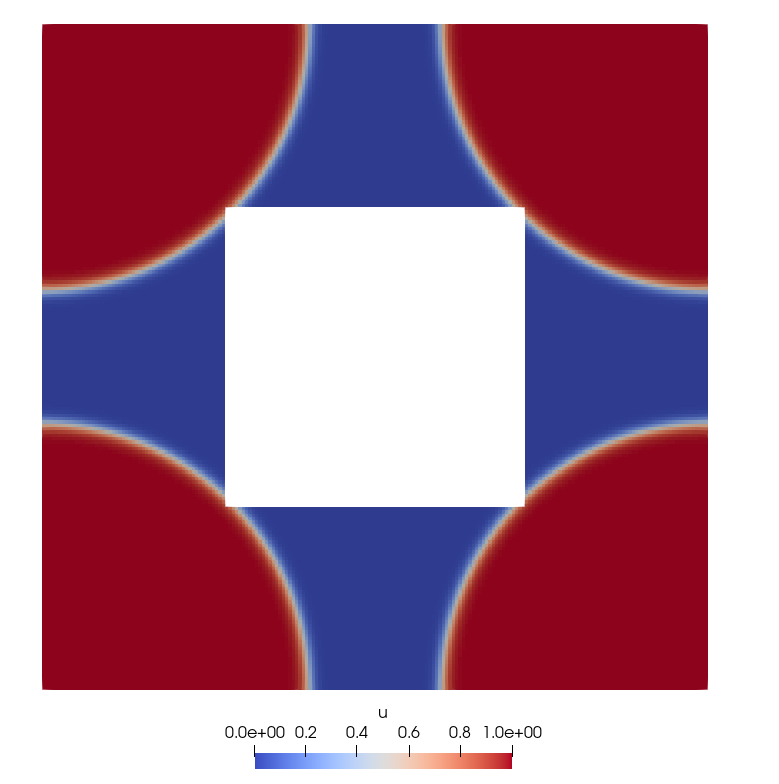}
\end{subfigure}
\begin{subfigure}{0.5\textwidth}
\includegraphics[width=\textwidth]{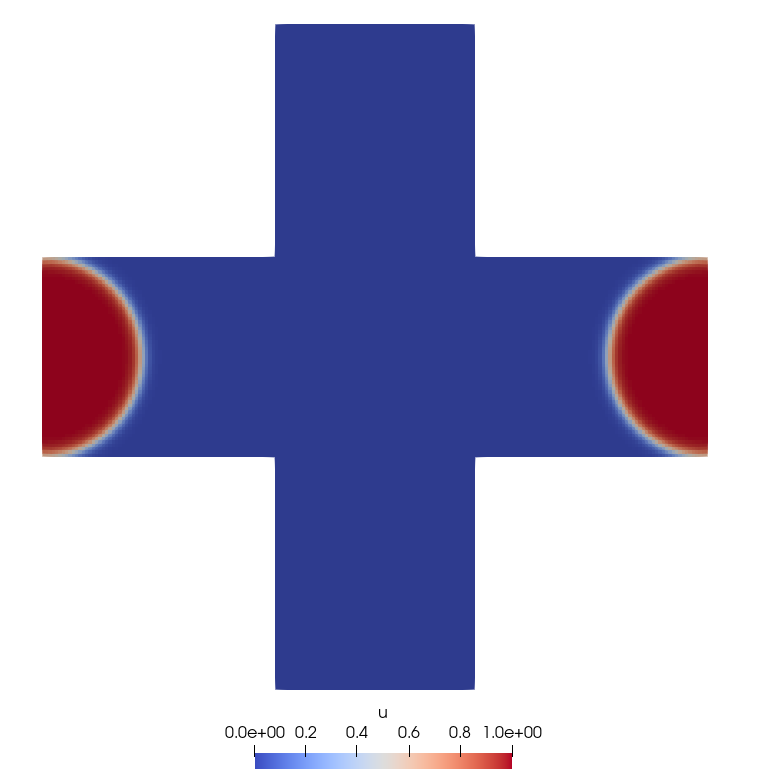}
\end{subfigure}
\caption{Investigated cell geometries "obstacle" and "cross". $\pf = 1$ (red) corresponds to fluid
    1, $\pf = 0$ (blue) to fluid 2.}
\label{fig:geometries}
\end{figure}

The two domains under consideration vary in geometry and porosity, see Figure~\ref{fig:geometries}.
The first setup, "obstacle", contains a square of side length $0.45$ placed in the center of the domain,
leaving a void fraction of $79.75\%$.
The second investigated geometry, denoted "cross", is a cross of channels with diameter $0.3$, resulting in a porosity
of $51\%$.

We consider four core scenarios of which the first three
investigate the effective mobility $\effPerm^{(k)}$
for different ratios of viscosity $\ratioVisc$ and density $\ratioDens$.
In the last setup the surface tension tensor $\effSurf^{(k)}$ is analyzed, with identical fluid
properties for both fluids.
The dimensionless numbers in the cell problems are chosen as 1.

\begin{figure}
\begin{subfigure}{0.5\textwidth}
\center
\includegraphics[width=\textwidth]{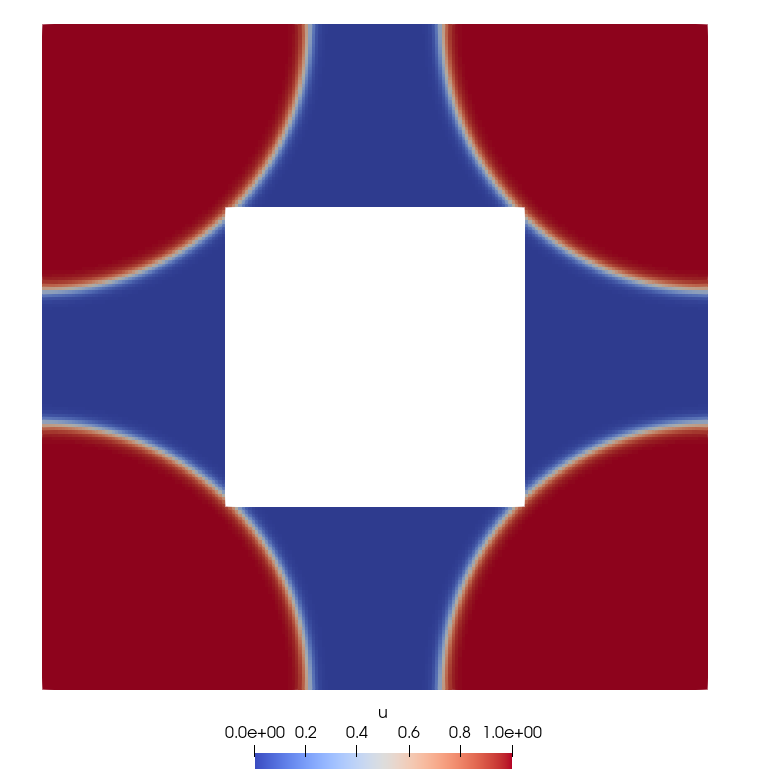}
\end{subfigure}
\begin{subfigure}{0.5\textwidth}
\center
\includegraphics[width=\textwidth]{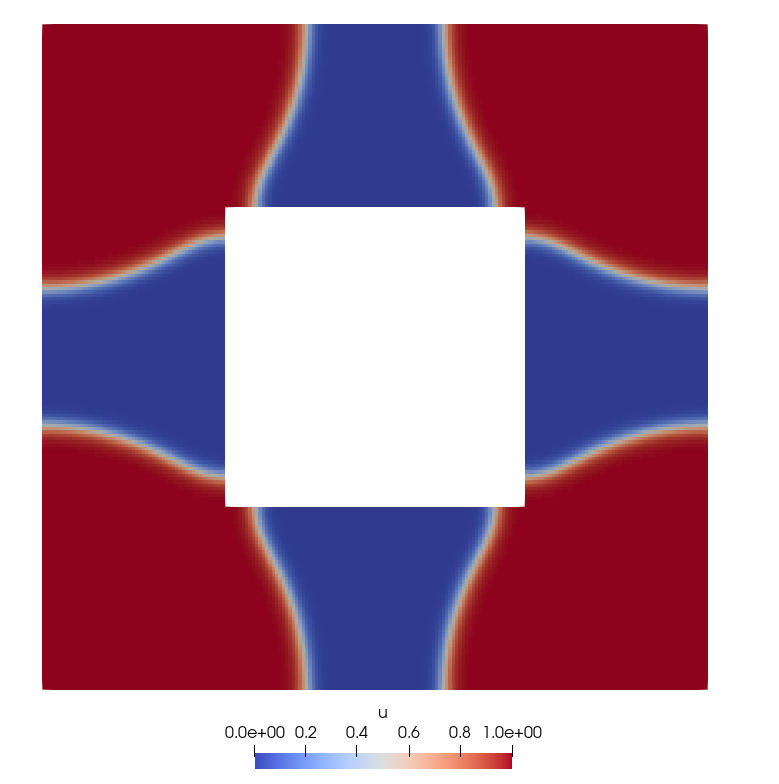}
\end{subfigure}
\caption{Preprocessing to enforce a neutral contact angle. Initial condition from
    \eqref{eq:pf_profile} (left), and the resulting phase field with desired contact angle (right).}
\label{fig:preprocessing}
\end{figure}

For all four cases we fix the center of a circle with varying radius $r$ and assign cells inside this
circle to phase 1, with the rest of the void space being filled with fluid 2.
This initial data already contains a diffuse interface according to the following radial function.
\begin{equation}
\label{eq:pf_profile}
\pf(r) = \frac{1}{1+\exp(5/\xi \cdot r)} ~.
\end{equation}
This initial data is evolved under the transient version of the Allen-Cahn equation
\eqref{eq:two_scale_cell_phasefield} without advection to enforce a
neutral contact angle (see Figure~\ref{fig:preprocessing}).
The resulting phase field is fed into the different cell problems for flow velocity components
\eqref{eq:two_scale_cell_permeability}, \eqref{eq:two_scale_cell_interface}, which are then solved.

In the "cross" cell geometry we observe recirculation patterns (Figure~\ref{fig:C10_paraview}), and e.g. for a horizontal pressure
gradient forcing flow from left to right the recirculation near the top and bottom of the domain
leads to flow from right to left at these edges of the periodic cell.
Depending on the fluid distribution one of the phases may contain primarily these flows in the
opposite direction to the main flow, leading to negative integrated velocities.
As the effective mobility should capture net flow through the cell, we exclude these
recirculations based on the computed flow field, avoiding negative mobilities.
The computed velocity field is loaded into ParaView
\cite{paraview}, computing streamlines from the boundaries and marking affected cells.
This process is not performed for the flow problem driven by surface tension
\eqref{eq:two_scale_cell_interface}.

After postprocessing we compute the weighted integrals over marked
cells to obtain the effective parameters.
The saturation and an approximation to the specific interfacial area are determined by integrals
over the entire domain.
\begin{equation}
S^{(1)} = \frac{1}{|\Po|} \int_\Po \pf_0 \dd\y \,, \quad
A = \frac{1}{|\Po|} \int_\Po \frac{4}{\xi} \pf_0 (1-\pf_0) \dd\y
\end{equation}

For the effective mobility $\effPerm^{(k)}$, the computed entries are divided by the absolute
permeability and multiplied with the viscosity, yielding a relative permeability.
The results are plotted over the saturation of phase 1.

\subsubsection{Equal fluid properties}
\begin{figure}
  \begin{subfigure}{0.5\textwidth}
  \includegraphics[width=\textwidth]{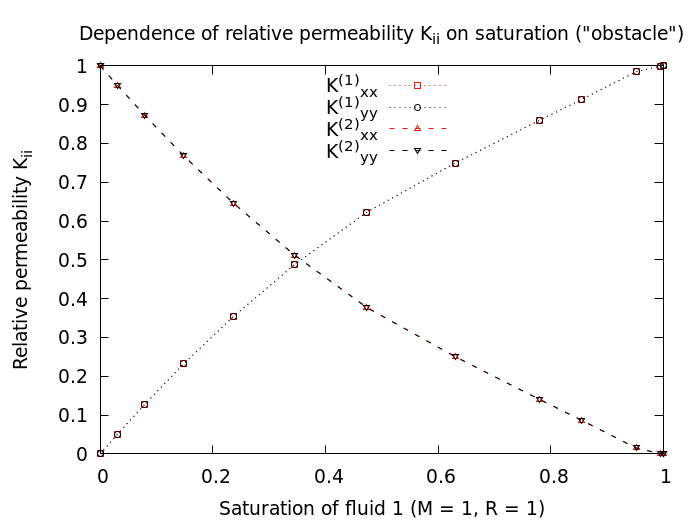}
  \end{subfigure}
  \begin{subfigure}{0.5\textwidth}
  \includegraphics[width=\textwidth]{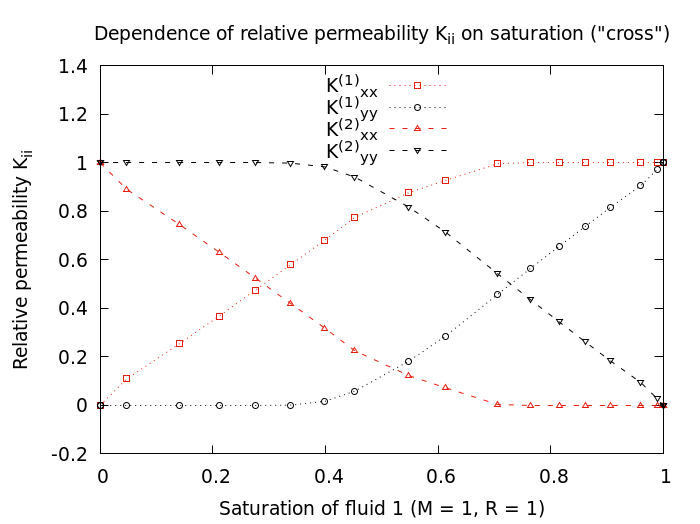}
  \end{subfigure}
  \caption{Relative-permeability-saturation relation for equal fluid properties.
  }
  \label{fig:A}
\end{figure}

In the first case we consider two fluids with equal viscosity and density
(\(\ratioVisc = \ratioDens = 1\)).
The cell problems for velocity contributions driven by the global pressure gradient treat the
two-phase system as single-phase flow.
The effective mobilities are obtained through integration
weighted by the phase-indicator $\pf^{(k)}$ and the resulting relative permeabilities always sum to
1, see Figure~\ref{fig:A}.

Due to the symmetries in the fluid distribution the relative permeability is approximately a diagonal matrix and
we show only the dominant diagonal entries.
For the isotropic setup in the "obstacle" geometry the entries are equal and the relative permability
reduces to a scalar value.

\begin{figure}
  \begin{subfigure}{0.5\textwidth}
  \includegraphics[width=\textwidth]{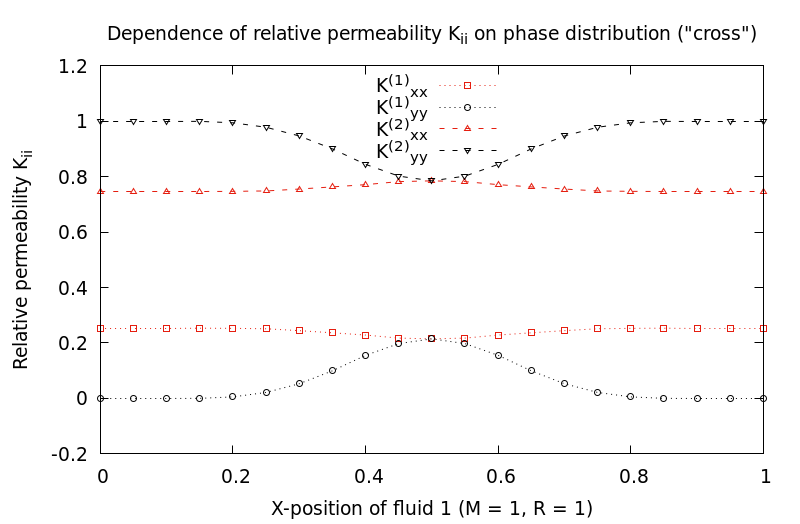}
  \end{subfigure}
  \begin{subfigure}{0.5\textwidth}
  \includegraphics[width=\textwidth]{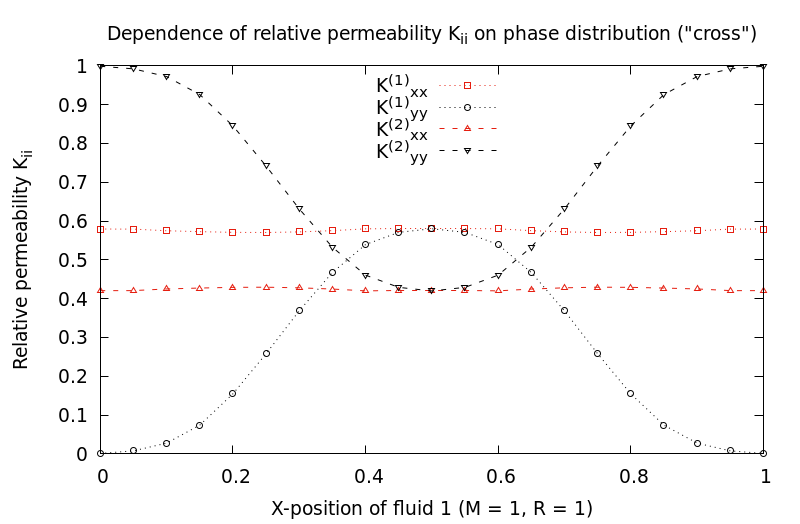}
  \end{subfigure}
  \caption{Impact of phase distribution on relative permeability for equal fluid properties with
      radius $r=0.15$ (left) and $r=0.3$ (right).}
  \label{fig:A_pos}
\end{figure}

Due to the choice of fluid distribution for every radius $r$ the domain corresponding to fluid 1
contains that for lower radii $r' < r$ and thus lower saturation.
Equivalently, for phase-field functions $\pf_1 < \pf_2$ the corresponding saturations fulfil
$S_1 < S_2$.
Under this condition the relative permeability is monotone with respect to the saturation.
For differing fluid distributions this monotonicity is not guaranteed and a higher saturation concentrated at areas
of lower velocities can result in a lower permeability.
Figure~\ref{fig:A_pos} shows changing relative permeabilities for fixed saturation and varying fluid
distributions obtained by changing the center of the initial circular subdomain for fluid 1.
This information about local phase distribution and its impact on macroscopic flow can only be
captured by solving the cell problems.

\subsubsection{Influence of viscosity differences}
Next we investigate the effect of the viscosity ratio by considering \(\ratioVisc = 2\), the
investigated cell problems now correspond to an incompressible fluid with varying viscosity.

\begin{figure}
  \begin{subfigure}{0.5\textwidth}
  \includegraphics[width=\textwidth]{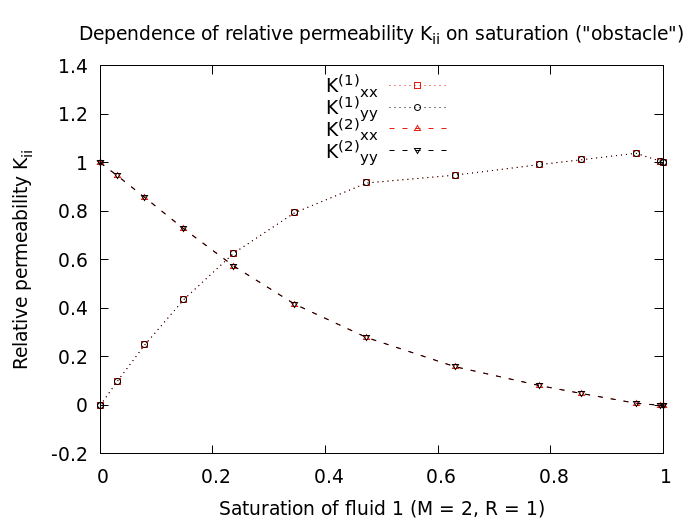}
  \end{subfigure}
  \begin{subfigure}{0.5\textwidth}
  \includegraphics[width=\textwidth]{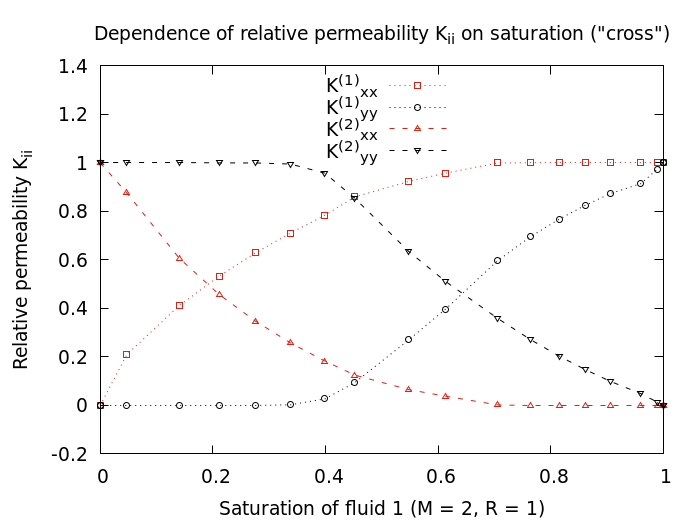}
  \end{subfigure}
  \caption{Relative-permeability-saturation relation for viscosity ratio \(\ratioVisc = 2\).}
  \label{fig:B}
\end{figure}

In stark contrast to typically used relative permeability curves we observe non-monotone behavior
for the case "obstacle" with
relative permeabilities of the more viscous fluid reaching values above 1 for saturations
above 0.8, up to $\effPerm^{(1)}_{xx} = 1.039$ at a saturation of 0.95 (see Figure~\ref{fig:B}). 
Due to the chosen fluid distribution the less viscous fluid 2 is located on the surface of the solid matrix
at low saturations (high saturation for phase 1).
This reduces the resistance exerted on the first fluid and results in higher velocities, especially
for higher viscosity ratios \(\ratioVisc\).
Such lubricating effects have been observed experimentally in different settings \cite{berg2008}.

\subsubsection{Influence of density differences}
\begin{figure}
  \begin{subfigure}{0.5\textwidth}
  \includegraphics[width=\textwidth]{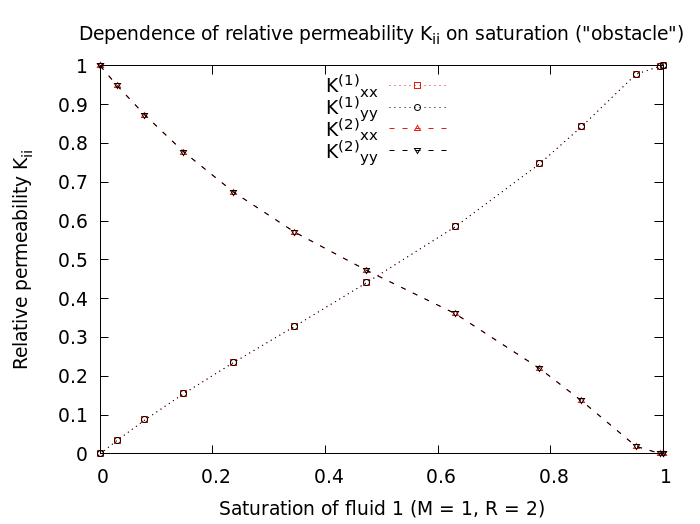}
  \end{subfigure}
  \begin{subfigure}{0.5\textwidth}
  \includegraphics[width=\textwidth]{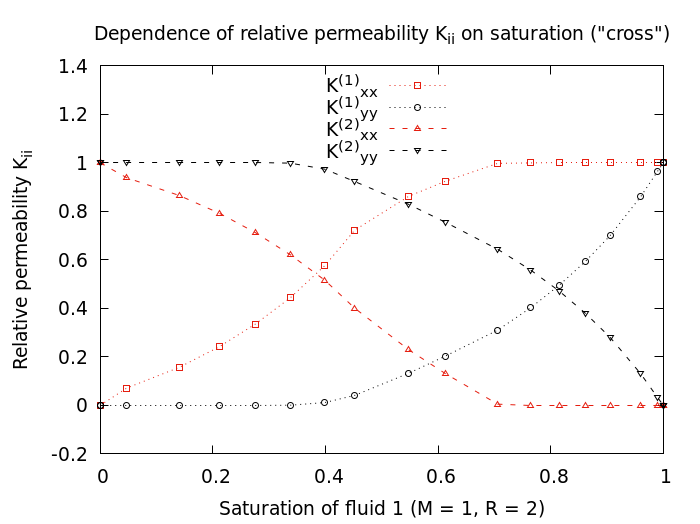}
  \end{subfigure}
  \caption{Relative-permeability-saturation relation for density ratio \(\ratioDens = 2\).}
  \label{fig:C2}
\end{figure}

\begin{figure}
  \begin{subfigure}{0.5\textwidth}
  \includegraphics[width=\textwidth]{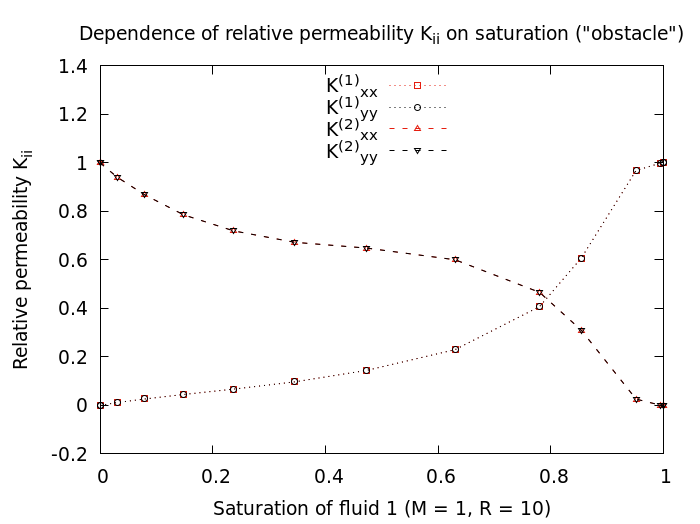}
  \end{subfigure}
  \begin{subfigure}{0.5\textwidth}
  \includegraphics[width=\textwidth]{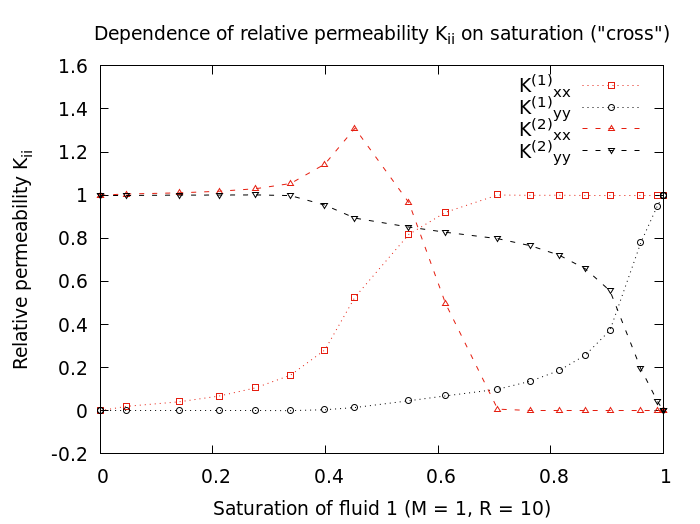}
  \end{subfigure}
  \caption{Relative-permeability-saturation relation for density ratio \(\ratioDens = 10\).}
  \label{fig:C10}
\end{figure}

\begin{figure}
  \includegraphics[width=\textwidth]{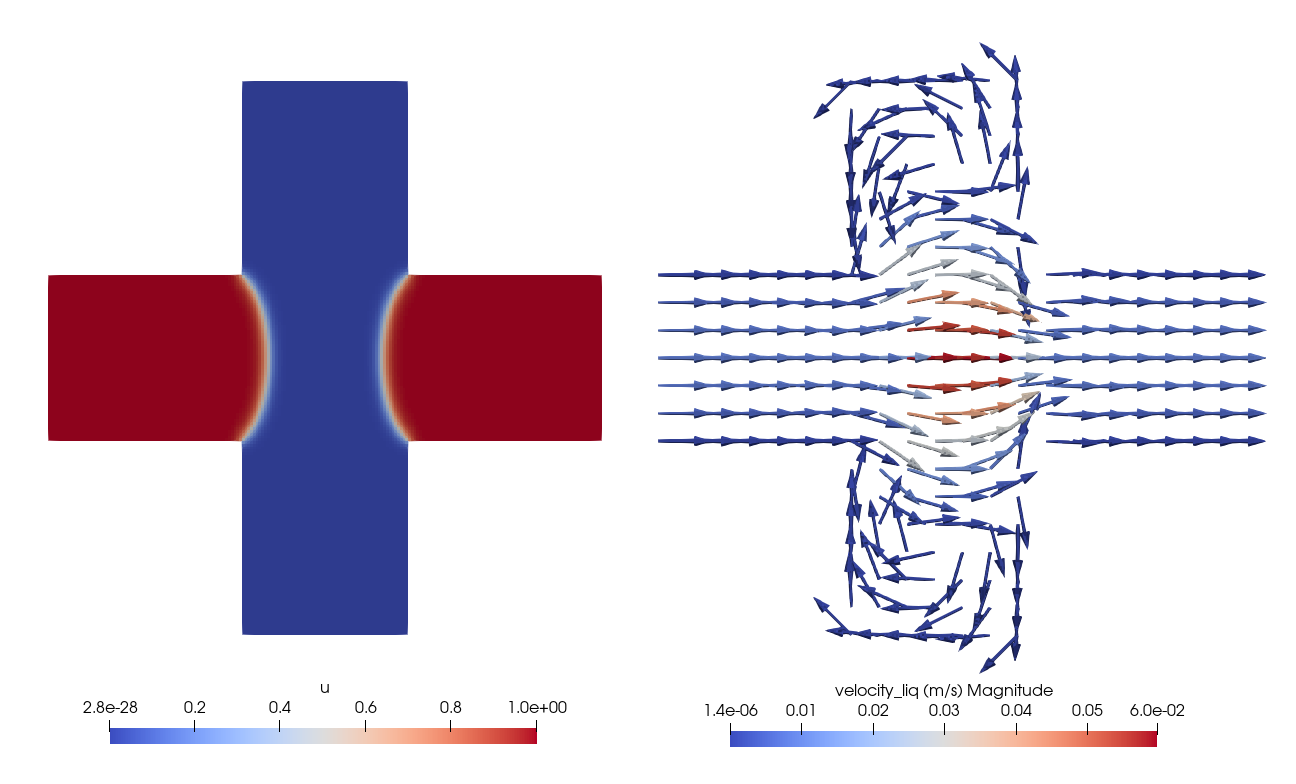}
  \caption{Phase distribution and velocity component for density ratio \(\ratioDens = 10\).}
  \label{fig:C10_paraview}
\end{figure}

\begin{figure}
  \begin{subfigure}{0.5\textwidth}
  \includegraphics[width=\textwidth]{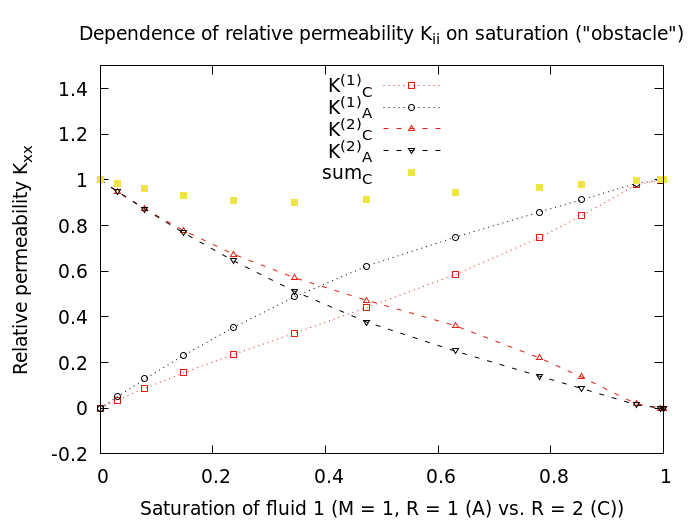}
  \end{subfigure}
  \begin{subfigure}{0.5\textwidth}
  \includegraphics[width=\textwidth]{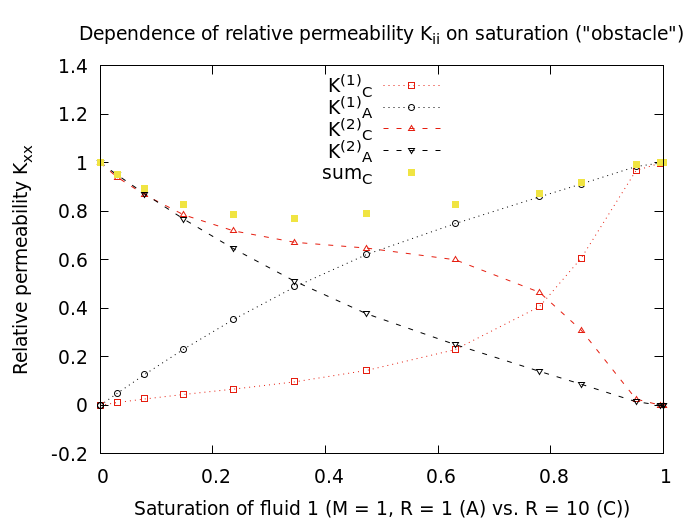}
  \end{subfigure}
  \caption{Comparison of relative permeability saturation relation for "obstacle" geometry.}
  \label{fig:C_comp_obstacle}
\end{figure}

In the case of non-trivial density ratio \(\ratioDens\), the cell problem
\eqref{eq:two_scale_cell_permeability} corresponds to Stokes equation for a quasi-compressible
fluid.
For ratio \(\ratioDens = 2\) (Figure~\ref{fig:C2})
no non-monotone behavior is observed but the curves are notably different to the one obtained for
the case of equal fluid properties.
For density ratio \(\ratioDens = 10\) (Figure~\ref{fig:C10}) non-monotonicity is observed for a
horizontal pressure gradient in the "cross" setup, with the velocity in the lighter fluid
increasing as the main flow path is filled with the more dense fluid.
The fluid distribution and computed flow field for highest permeability at a saturation of 0.45
are shown in Figure~\ref{fig:C10_paraview}.
In Figure~\ref{fig:C_comp_obstacle} the results for the "obstacle"
geometry are compared to those for equal fluid properties, which highlights the fact that the sum of
the relative permeabilities is less than one.


\subsubsection{Surface tension tensor}
\begin{figure}
  \begin{subfigure}{0.5\textwidth}
  \includegraphics[width=\textwidth]{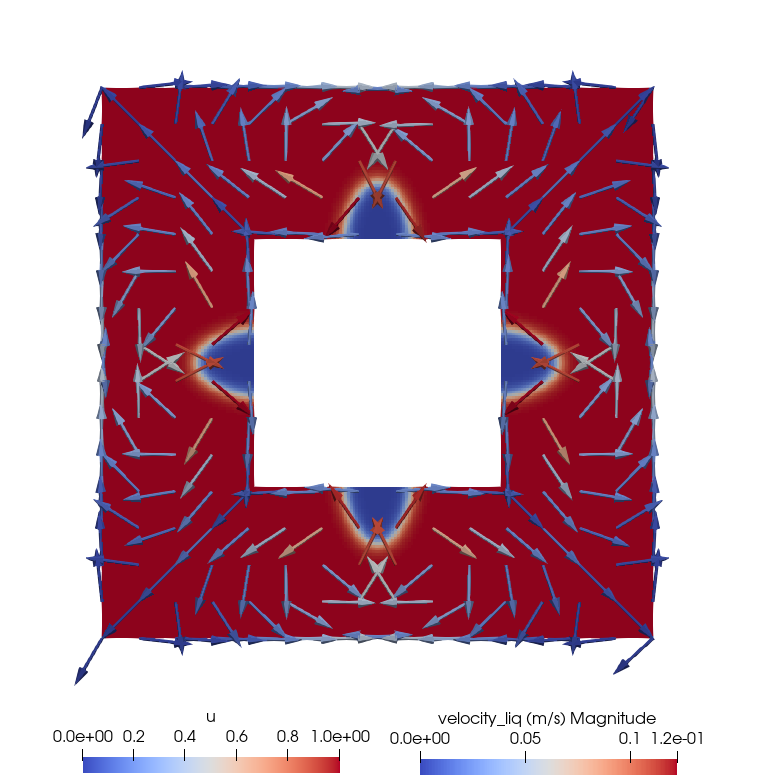}
  \end{subfigure}
  \begin{subfigure}{0.5\textwidth}
  \includegraphics[width=\textwidth]{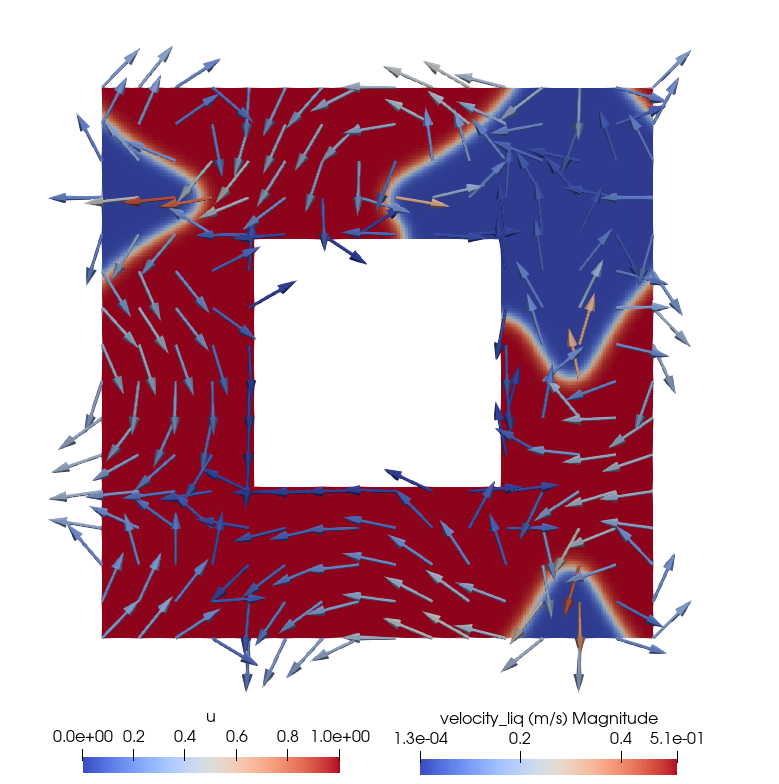}
  \end{subfigure}
  \caption{Symmetric and asymmetric case with center at (0.35, 0.3).}
  \label{fig:asymm_both}
\end{figure}

For the isotropic geometry and fluid distribution considered above the velocities driven by
surface tension cancel out in the integral and the computed effective parameter is equal to 0
(Figure~\ref{fig:asymm_both}, left).
For asymmetric fluid distributions (Figure~\ref{fig:asymm_both}, right) a
non-zero contribution is obtained, but no visible trends are observed as the saturation changes.
The redistribution of fluids driven by surface tension forces leads to a net flow for the two
phases, which corresponds to $\effSurf^{(k)}$ being non-zero.
The size and direction of this net flow and hence the size and sign of the effective parameter
$\effSurf^{(k)}$ depends highly on the fluid distribution.
Such effects can only be accounted for by solving the cell problem
\eqref{eq:two_scale_cell_interface}.

\subsection{Discussion}
\label{sec:discussion}
For equal fluid properties and fixed fluid distributions we obtain monotone
relative permeability saturation curves comparable to commonly used relations such as Brooks-Corey
\cite{brooks1966} or van Genuchten \cite{vangenuchten1980}.
However, even in this simple case the relative permeability depends not only on saturation but also
on the distribution of the fluids.
As observed in Section~\ref{sec:numerics_cell}, even in this simple case no relation can be given without accounting for the local
distribution of the fluids.

For differing viscosities the model is able to reproduce non-monotone relative permeabilities with
values above 1 as discussed in \cite{berg2008}. 
When the less viscous fluid coats the solid surface it can reduce the resistance exerted on the
more viscous fluid, resulting in higher velocities.
Due to the different fluid distribution such a coating is not observed to the
same extent for the "cross" setup, leading to monotone relations for the relative permeability.
These effects depend strongly on the local phase distribution and cannot be captured solely by the
saturation without additional assumptions.

For moderate density differences the combined permeability is reduced and for higher density ratios
non-monotone behavior can be observed.
We remark that for very high differences in fluid properties, this should be considered in the derivation of the cell
problems and the two-scale problem, which would lead to a different appearance of the
effective equations.
To obtain equations similar to the often used two-phase Darcy's law, we use the assumptions
presented in Section~\ref{sec:homogenization}.

The derived two-scale model contains an additional effective parameter which accounts for the
influence of surface tension between the fluid phases.
This is not a common part of the extended Darcy's law and there is no counterpart to compare it to.
A similar cell problem and resulting effective parameter has been investigated by
\cite{sharmin2022}, incorporating the influence of solute-dependent surface tension.
Here only constant surface tension is considered.
For symmetric distributions the effective parameter capturing surface tension effects disappears as
expected.
For anisotropic phase distributions the flow driven by surface tension does not cancel out and
the associated cell problem is able to capture significant
contributions to effective flow behavior.
The size and direction of flow resulting from surface tension and thus the magnitude and sign of the
effective parameter depend strongly on the local distribution of fluids.

The numeric investigation of the effective parameters highlights how different fluid distributions
with equal saturation can result in very different net flow and effective behavior.
Increased relative permeability due to viscosity differences and significant surface tension effects further
emphasize the importance of resolving the local fluid morphology.

\section{Conclusion}
\label{sec:conclusion}
We derived a two-scale model for two-phase flow in porous media using homogenization and investigated the dependence of
effective parameters on the fluid distribution in the pore-scale cell problems. 

The model consists of macro-scale equations similar to the extended Darcy's law with effective
parameters computed from solutions to cell problems defined at every macroscopic point.
One of the effective parameters can be understood as an effective mobility, while the second
effective parameter represents effects of interfacial tension on the effective flow.
On the pore-scale the phase distribution is captured by a phase field, using an advective Allen-Cahn formulation.
For the velocity field we use a stationary Stokes equation with phase-dependent fluid properties and
an additional source term incorporating surface tension forces.

The two-scale structure allows to investigate numerically the effects of pore-scale information on
effective behavior.
We construct local fluid distributions for the pore scale and solve the cell problems associated
with pressure driven velocity contributions.

The similarity to two-phase Darcy's law motivates a comparison of the dependence of relative permeabilities on
saturation. By selecting isotropic pore geometries, relative permeability tensors can be extracted
from the effective parameters \(\effPerm^{(k)}\) in our model and for isotropic fluid distributions
these simplify to diagonal matrices or scalar values.
As one of the core qualities of relative-permability-saturation-curves we investigate the conditions
under which a monotone relationship is obtained.
If a common fluid morphology is maintained for different saturations, then for equal fluid
properties monotone curves, similar to commonly used functions like van Genuchten and Brooks-Corey,
are obtained.
The calculated curves however depend strongly on the fluid distribution and change with the pore geometry,
highlighting the need to account for pore-scale information.
For different viscosities the model can produce non-monotone relationships and relative
permeabilities greater than 1. This can be related to a lubricating effect of the less viscous fluid
coating the surface of the solid matrix, an effect observed experimentally.
For moderate density ratios the effective total permeability decreases and no
non-monotone behavior was observed.

We also performed numeric simulations for the cell problem corresponding to the effective parameter
\(\effSurf^{(k)}\), capturing effects of surface tension forces.
For isotropic geometry and fluid morphologies these forces result in no net flow and a vanishing
effective parameter.
For asymmetric fluid distributions we obtain significant net flows captured by \(\effSurf^{(k)}\).

Our two-scale model for two-phase flow captures pore-scale fluid-fluid interactions related to
surface tension forces as well as differences in viscosity and density.
Through numeric investigation of the effective parameters we highlight the importance of local fluid
distribution for accurately capturing effective flow behavior.

\paragraph{Acknowledgments}
Funded by Deutsche Forschungsgemeinschaft (DFG, German Research Foundation) under Germany's
Excellence Strategy - EXC 2075 - 390740016. We acknowledge the support by the Stuttgart Center for
Simulation Science (SimTech).
We thank the Deutsche Forschungsgemeinschaft (DFG, German Research Foundation) for supporting this work by funding SFB 1313, Project Number 327154368.

The authors would like to thank Kundan Kumar (University of Bergen) for useful discussions on
monotonicity of relative permeability saturation curves.

\appendix

\section{Sharp-interface limit}\label{app:si}

In the limit $\xi \to 0$ the phase-field model presented in Section~\ref{sec:pf_model} recovers the
sharp-interface model from Section~\ref{sec:si_model}.
Following the approach of matched asymptotic expansions \cite{caginalp1988} we consider the behavior
at the transition zone and in
the bulk phases separately, connecting them with matching conditions.
Let $\refL$ be a reference length and $\xxi = \xi/\refL$ the non-dimensional interface width.

We assume asymptotic expansions in $\xxi$ that capture the behavior away from the interface (outer
expansions) for unknowns $\psi \in \{ \pf, p, \vv \}$
\begin{equation}
\label{eq:outer_expansion}
\psi^\out(t, \x) = \sum_{k=0}^\infty \xxi^k \psi^\out_k (t, \x) ~.
\end{equation}
For inner expansions at the diffuse interface we consider local curvilinear coordinates.
Let $\Interface_\xxi (t) = \{ \y_\xxi \in \Po \mid \pf(t, \y_\xxi) = 1/2 \}$ be the reconstructed interface.
With a parametrization $\s$ along $\Interface_\xxi(t)$ and normal $\n_\xxi(t, \s)$ pointing into fluid
2, we define a signed distance $r$ for points near the interface, such that
\begin{equation}
\label{eq:curvilinear_coordinates}
\x = \y_\xxi(t, \s) + r \n_\xxi(t, \s) ~.
\end{equation}
Note than $\partial_t r = -v_n$, with $v_n$ the velocity of the interface in direction of $\n_\xxi$
corresponding to the sharp-interface velocity $V_\InterfaceFluids$.
Furthermore it can be shown (see \cite{caginalp1988}) that for mean and Gaussian interface curvatures
$\kappa$ and $\Pi$
\begin{equation*}
|\nabla r| = 1,\quad \nabla r \cdot \nabla s_i = 0,
\quad \nabla^2 r = \frac{\kappa + 2\Pi r}{1+\kappa r +\Pi r^2} ~.
\end{equation*}
The outer expansions define the interface
$\Interface^\out_0(t) = \{ \x \in \Po \mid \pf^\out_0(t, \x) = 1/2\}$ with normal vector $\n_0$,
interfacial velocity $v_{n,0}$ and mean curvature $\kappa_0$.
The point $\y_\xxi$ can be expressed through expansions
\begin{equation*}
\y_\xxi = \sum_{k=0}^\infty \xxi^k \y_k
\end{equation*}
with $\y_0 \in \Interface^\out_0(t)$, and similarly $\n_\xxi = \n_0 + \xxi \n_1 + \Ox{2}$.
Defining $z = r/\xxi$, we consider the inner expansions in the variables $z$ and $\s$
\begin{equation}
\label{eq:inner_expansion}
\psi^\inn(t, \x) = \sum_{k=0}^\infty \xxi^k \psi^\inn_k(t, z, \s) ~.
\end{equation}
The derivatives are rewritten accordingly, yielding \cite{caginalp1988}
\begin{equation}
\begin{aligned}
\label{eq:inner_derivatives}
\partial_t \psi &= -\xxi^{-1} v_{n,0} \partial_z \psi^\inn
+ (\partial_t + \partial_t \s \cdot \nabla_\s) \psi^\inn
+ \Ox{}
~, \\
\nabla_\x \psi &= \xxi^{-1} \partial_z \psi^\inn \n_0 + \nabla_\Interface \psi^\inn + \Ox{}
~, \\
\nabla_\x \cdot \ve{\psi} &= \xxi^{-1} \partial_z \ve{\psi}^\inn \cdot \n_0
+ \nabla_\Interface \cdot \ve{\psi}^\inn + \Ox{}
~, \\
\nabla_\x^2 \psi &= \xxi^{-2} \partial_{zz} \psi^\inn + \xxi^{-1} \kappa_0 \partial_z \psi^\inn
+ \mathcal{O}(1)
~,
\end{aligned}
\end{equation}
using the expansions $v_n = v_{n,0} + \Ox{}$ and $\kappa = \kappa_0 + \Ox{}$.
For fixed $t$ and $\s$, let $\y_{1/2\pm}$ denote the limit for $r\searrow 0$ and $r\nearrow 0$
respectively of a point $\y$ in curvilinear coordinates \eqref{eq:curvilinear_coordinates}.
The values of outer expansions at $\y_{1/2\pm}$ are paired with the limit values of inner expansions
for $z \to \pm\infty$ in the following matching conditions.
\begin{equation}
\begin{aligned}
\label{eq:matching_conditions}
\lim_{z\to \pm \infty} \psi^\inn_0 (t, z, \s) &= \psi^\out_0 (t, \y_{1/2\pm}) ~, \\
\lim_{z\to\pm\infty} \partial_z \psi^\inn_0(t, z, \s) &= 0 ~, \\
\lim_{z\to\pm\infty} (\psi^\inn_0(t, z, \s) - (z + y_1)\nabla \psi^\out_0(t, \y_{1/2\pm}) \cdot \n_0)
&= \psi^\out_1(t, \y_{1/2\pm}) ~, \\
\lim_{z\to\pm\infty} \partial_z \psi^\inn_1(t, z, \s) &= \nabla \psi^\out_0(t, \y_{1/2\pm}) \cdot \n_0 ~.
\end{aligned}
\end{equation}

\subsection{Outer expansions}
Due to the polynomial structure of $P'$ it can be expanded as
\begin{equation}
P'(\pf_0 + \xxi \pf_1 + \Ox{2}) = P'(\pf_0) + \xxi \pf_1 P''(\pf_0) + \Ox{2} ~.
\end{equation}
Substituting the outer expansions \eqref{eq:outer_expansion} into the phase-field equation
\eqref{eq:pore_phasefield} yields
\begin{subequations}
\begin{alignat}{1}
0 &= \refL^2 \xxi^2 \frac{\partial}{\partial t} (\pf^\out_0 + \Ox{1})
    + \refL^2 \xxi^2 \nabla \cdot ((\vv_0 + \Ox{1})(\pf^\out_0 + \Ox{1}))
\notag\\&
    - \mob \refL^3 \xxi^2 \nabla^2 (\pf^\out_0 + \Ox{1})
    + \mob \refL \xxi P'(\pf^\out_0 + \xxi \pf^\out_1 + \Ox{2})
\\&=
    \big[ \mob \refL \xxi P'(\pf^\out_0) \big]
    + \Ox{2}
~.
\end{alignat}
\end{subequations}

The leading order equation $P'(\pf^\out_0) = 0$ has solutions $\pf^\out_0 \in \{ 0, 1/2, 1 \}$,
of which 0 and 1 are stable minimizers of $P$.
This allows a decomposition of the unified $\Po$ into bulk domains $\Om^\out_1$ and $\Om^\out_2$ for
fluid 1 and 2, respectively, through
\begin{equation}
\Om^\out_1 (t) := \{ \x \in \Po \mid \pf^\out_0 (t, \x) = 1 \} ~,
    \qquad
\Om^\out_2 (t) := \{ \x \in \Po \mid \pf^\out_0 (t, \x) = 0 \} ~.
\end{equation}
We consider the remaining outer expansions for $\Om^\out_i$.

For the mass balance equation \eqref{eq:diff_mass} we obtain by using the linear dependence of
density on the phase-field variable
\begin{subequations}
\begin{alignat}{1}
0 =& \frac{\partial}{\partial t} (\rho ( \pf^\out_0 ) + \Ox{})
    + \nabla \cdot ( (\rho (\pf^\out_0) + \Ox{}) (\vv^\out_0 + \Ox{}) )
\\
  =& \frac{\partial \rho_i}{\partial t}
    + \nabla \cdot ( \rho_i \vv^\out_0)
    + \Ox{}
    ~.
\end{alignat}
\end{subequations}
Denoting phase velocities $\vv_i = {\vv^\out_0}\vert_{\Om^\out_i}$ we recover the individual mass
conservation equations and, due to constant phase densities $\rho_i$, the divergence-free condition
\eqref{eq:si_mass}.
The momentum balance \eqref{eq:diff_momentum} yields, using additionally the linear relation
for viscosity,
\begin{subequations}
\begin{alignat}{1}
0 =& \frac{\partial}{\partial t} ((\rho ( \pf^\out_0 ) + \Ox{}) (\vv^\out_0 + \Ox{}) )
\notag\\&
    + \nabla \cdot ( (\rho (\pf^\out_0) + \Ox{}) (\vv^\out_0 + \Ox{}) \otimes (\vv^\out_0 + \Ox{}) )
\notag\\&
    + \nabla (p^\out_0 + \Ox{})
\notag\\&
    - \nabla \cdot \bigg( (\mu (\pf^\out_0) + \Ox{})
    \bigg( \nabla (\vv^\out_0 + \Ox{}) + (\nabla (\vv^\out_0 + \Ox{}) )^T
\notag\\&
            -\frac{2}{3} (\nabla \cdot (\vv^\out_0 + \Ox{}) ) \tens{I} \bigg) \bigg)
\notag\\&
    - (\rho (\pf^\out_0) + \Ox{}) \ve{g}
    + \Ox{}
\\
  =& \frac{\partial}{\partial t} (\rho_i \vv_i )
    + \nabla \cdot ( \rho_i \vv_i \otimes \vv_i )
    + \nabla p^\out_0
\notag\\&
    - \nabla \cdot \bigg( \mu_i
    \bigg( \nabla \vv_i + (\nabla \vv_i)^T
            -\frac{2}{3} (\nabla \cdot \vv_i) \tens{I} \bigg) \bigg)
    - \rho_i \ve{g}
    + \Ox{}
\end{alignat}
\end{subequations}
With phase pressures $p_i$ given by restricting $p^\out_0$ to the bulk domain $\Om^\out_i$ we
recover the momentum balance \eqref{eq:si_momentum}.

Inserting the outer expansions into the boundary conditions yields the conditions for the
fluid-solid interfaces of the sharp-interface model.
From \eqref{eq:diff_bc_slip} we obtain
\begin{subequations}
\begin{alignat}{1}
0 &= (\vv^\out_0 + \Ox{}) + \slip \nabla ((\vv^\out_0 + \Ox{}) \cdot (\ve{t}_0 + \Ox{}))
    \cdot (\n_0 + \Ox{}) (\ve{t}_0 + \Ox{})
\\
&= \vv_i + \slip \nabla (\vv_i \cdot \ve{t}_0) \cdot \n_0 \,\ve{t}_0
~,
\end{alignat}
\end{subequations}
corresponding to the slip condition \eqref{eq:si_bc_slip}.

\subsection{Inner expansions}
To obtain the boundary conditions at the fluid-fluid interface we insert the inner expansions
\eqref{eq:inner_expansion}, rewriting derivatives according to \eqref{eq:inner_derivatives}, and use
the matching conditions \eqref{eq:matching_conditions}.

\paragraph{Phase-field equation}
From the phase-field equation \eqref{eq:pore_phasefield} we obtain due to the polynomial form of
$P$
\begin{subequations}
\begin{alignat}{1}
0 =& \refL^2 \xxi^2 (-\xxi^{-1} v_{n,0} \partial_z (\pf^\inn_0 + \Ox{}))
    + \refL^2 \xxi^2 (-\xxi^{-1} \partial_z ((\vv^\inn_0 + \Ox{}) (\pf^\inn_0 + \Ox{}))
        \cdot \n_0)
\notag\\&
    - \mob \refL^3 \xxi^3 (\xxi^{-2} \partial_{zz} (\pf^\inn_0 + \xxi \pf^\inn_1 + \Ox{2})
        + \xxi^{-1} \kappa_0 \partial_z (\pf^\inn_0 + \Ox{}))
\notag\\&
    + \mob \refL \xxi (P'(\pf^\inn_0) + \pf^\inn_1 P''(\pf^\inn_0) + \Ox{2})
\\=&
    \xxi \big[ -\mob \refL^2 \partial_{zz} \pf^\inn_0 + \mob P'(\pf^\inn_0)
    -\refL^2 v_{n,0} \partial_z \pf^\inn_0
    + \refL^2 \partial_z (\vv^\inn_0 \pf^\inn_0) \cdot \n_0
    \big]
    + \Ox{2}
~.
\end{alignat}
\end{subequations}
Multiplying the leading order terms with $\partial_z \pf^\inn_0$ and integrating with respect to $z$ yields a
condition for the interface velocity
\begin{subequations}
\begin{alignat}{1}
0 =& \int_{z=-\infty}^{\infty}
\mob (-\refL^3 \partial_{zz} \pf^\inn_0 + \refL P'(\pf^\inn_0)) \partial_z \pf^\inn_0
+ (-v_{n,0} + \vv^\inn_0 \cdot \n_0) \refL^2 (\partial_z \pf^\inn_0)^2 \dd z
\\=& \mob \refL^2 \int_{z=-\infty}^{\infty}
- \frac{\refL}{2} \partial_z (\partial_z \pf^\inn_0)^2
+ \partial_z P(\pf^\inn_0) \dd z
+ (-v_{n,0} + \vv^\inn_0 \cdot \n_0) \refL^2
\int_{z=-\infty}^{\infty} (\partial_z \pf^\inn_0)^2 \dd z
\\=& \mob \refL^2 \big[ -\frac{L}{2} (0 - 0) + (0 - 0) \big] +
(-v_{n,0} + \vv^\inn_0 \cdot \n_0) \frac{2\refL}{3}
~,
\end{alignat}
\end{subequations}
where the terms in the first integral are equal to 0 due to the matching conditions.
Taking the limit in $z$, we obtain
\begin{equation}
v_{n,0} = \vv^\out_0(t, \y_{1/2\pm}) \cdot \n_0 ~.
\end{equation}

With $\vv^\inn_0 \cdot \n_0$ independent of $z$, the leading order terms reduce to
\begin{equation}
0 = -\mob \refL^2 \partial_{zz} \pf^\inn_0 + \mob P'(\pf^\inn_0) ~.
\end{equation}
Multiplying the remaining terms with $\partial_z \pf^\inn_0$ and integrating with respect to
$z$ yields, applying the chain rule and matching conditions with $\pf^\inn_0 (t, \y_{1/2-}) = 0$
and $P(0) = 0$,
\begin{subequations}
\begin{alignat}{1}
0 =& \int_{z'=-\infty}^{z} -\refL^2 \partial_z \pf^\inn_0 \partial_{zz} \pf^\inn_0
    + \partial_z \pf^\inn_0 P'(\pf^\inn_0) \dd z
\\=& \int_{z'=-\infty}^{z} -\refL^2 \frac{1}{2} \partial_z ((\partial_z \pf^\inn_0)^2)
    + \partial_z P(\pf^\inn_0) \dd z
\\=& -\refL^2 (\partial_z \pf^\inn_0(t, z, \ve{s}))^2 - 0 + 2 P(\pf^\inn_0(t, z, \ve{s})) - 0
~.
\end{alignat}
\end{subequations}
By definition of the curvilinear coordinates we have $\partial_z \pf^\inn_0 \ge 0$ and
\begin{equation}
\partial_z \pf^\inn_0 = \refL^{-1} \sqrt{ 2 P(\pf^\inn_0) } ~.
\end{equation}
Using $\pf^\out(t, 0, \ve{s}) = 1/2$ one can obtain a profile for the phase-field, namely
\begin{equation}
\pf^\inn_0(t, z, \ve{s}) = \pf^\inn_0(z) = \frac{1}{1+\exp{4z/\refL}}
= \frac{1}{2} \left( 1 + \tanh \left(\frac{2z}{\refL}\right)\right)
    ~.
\end{equation}

\paragraph{Mass balance}
Inserting the inner expansions into the mass balance equation \eqref{eq:diff_mass} yields
\begin{subequations}
\begin{alignat}{1}
0 =& (-\xxi^{-1} v_{n,0} \partial_z + \mathcal{O}(1)) (\rho(\pf^\inn_0) + \Ox{})
\notag\\&
    + \xxi^{-1} \partial_z ((\vv^\inn_0 + \Ox{})(\rho(\pf^\inn_0) + \Ox{})) \cdot \n_0 + \mathcal{O}(1)
\\=& \xxi^{-1} \big[ -v_{n,0} \partial_z \rho(\pf^\inn_0)
    + \partial_z (\vv^\inn_0 \rho(\pf^\inn_0)) \cdot \n_0 \big] + \mathcal{O}(1)
~.
\end{alignat}
\end{subequations}
Integrating with respect to $z$ and applying the matching conditions we obtain
\begin{subequations}
\begin{alignat}{1}
0 =& \int_{z=-\infty}^{\infty} -v_{n,0} \partial_z \rho(\pf^\inn_0)
    + \partial_z (\vv^\inn_0 \rho(\pf^\inn_0)) \dd z
\\=& -v_{n_0} (\rho(1) - \rho(0))
    + (\vv^\out_0(t, \y_{1/2+}) \rho(1) - \vv^\out_0(t, \y_{1/2-}) \rho(0))
\\=& -v_{n,0} (\rho_1 - \rho_2)
    + (\vv_1 \rho_1 - \vv_2 \rho_2)
~,
\end{alignat}
\end{subequations}
which is fulfilled for $v_1 = v_2 = v_{n,0}$ at the interface.

\paragraph{Momentum balance}
Inserting the inner expansions into the momentum balance equation \eqref{eq:diff_momentum} yields
\begin{subequations}
\begin{alignat}{1}
0 =& -\xxi^{-1} v_{n,0} \partial_z ((\rho(\pf^\inn_0) + \Ox{}) (\vv^\inn_0 + \Ox{}))
    + \mathcal{O}(1)
\notag\\&
    + \xxi^{-1} \partial_z ((\rho(\pf^\inn_0) + \Ox{}) (\vv^\inn_0 + \Ox{}) \otimes (\vv^\inn_0 +
                \Ox{})) \n_0
    + \mathcal{O}(1)
\notag\\&
    + \xxi^{-1} \partial_z (p^\inn_0 + \Ox{}) \n_0 + \mathcal{O}(1)
\notag\\&
    - \xxi^{-1} \partial_z \bigg(
        (\mu(\pf^\inn_0) + \xxi(M-1)\pf^\inn_1 + \Ox{2}) \bigg(
            \xxi^{-1} \partial_z (\vv^\inn_0 + \xxi \vv^\inn_1 + \Ox{2}) \otimes \n_0
\notag\\&
            + \nabla_\Interface (\vv^\inn_0 + \Ox{})
            + \xxi^{-1} (\partial_z (\vv^\inn_0 + \xxi \vv^\inn_1 + \Ox{2}) \otimes \n_0)^T
            + (\nabla_\Interface (\vv^\inn_0 + \Ox{}))^T
\notag\\&
            - \xxi^{-1} \frac{2}{3} \partial_z (\vv^\inn_0 + \xxi \vv^\inn_1 + \Ox{2})\cdot \n_0 \tens{I}
            - \frac{2}{3} \nabla_\Interface \cdot (\vv^\inn_0 + \Ox{}) \tens{I}
             + \Ox{}
    \bigg) \bigg) \n_0
\notag\\&
    + \nabla_\Interface \cdot \bigg( (\mu(\pf^\inn_0) +\Ox{}) \xxi^{-1} \bigg(
            \partial_z (\vv^\inn_0 + \Ox{}) \otimes \n_0
            + (\partial_z (\vv^\inn_0 + \Ox{}) \otimes \n_0)^T
\notag\\&
            - \frac{2}{3} \partial_z (\vv^\inn_0 + \Ox{}) \cdot \n_0 \tens{I}
            + \Ox{} \bigg) \bigg)
    - (\rho(\pf^\inn_0) + \Ox{}) \ve{g}
\notag\\&
    - \frac{3 \tension \refL}{2} \xxi \Big[
        \xxi^{-1} \partial_z \Big(
            \xxi^{-1} \partial_z (\pf^\inn_0 + \xxi \pf^\inn_1 + \Ox{2}) \n_0
             + \nabla_\Interface (\pf^\inn_0 + \Ox{})
\notag\\&
            \otimes
            (\xxi^{-1} \partial_z (\pf^\inn_0 + \xxi \pf^\inn_1 + \Ox{2}) \n_0
             + \nabla_\Interface (\pf^\inn_0 + \Ox{}) ) \Big) \n_0
\notag\\&
        + \nabla_\Interface \cdot \Big(
            \xxi^{-1} \partial_z (\pf^\inn_0 + \Ox{}) \n_0 \otimes
            \xxi^{-1} \partial_z (\pf^\inn_0 + \Ox{}) \n_0
        \Big)
    + \mathcal{O}(1) \Big]
\notag\\&
+ \mathcal{O}(1)
\intertext{with leading order}
=& \xxi^{-2} \bigg[
    -\partial_z ( \mu(\pf^\inn_0) ( \partial_z \vv^\inn_0 \otimes \n_0
                + (\partial_z \vv^\inn_0 \otimes \n_0)^T
                - \frac{2}{3} \partial_z \vv^\inn_0 \cdot \n_0 \tens{I} ) ) \n_0
\notag\\&
    - \frac{3 \tension \refL}{2} \partial_z (
            \partial_z \pf^\inn_0 \n_0 \otimes
            \partial_z \pf^\inn_0 \n_0) \n_0
    \bigg]
\intertext{and second order}
&
+\xxi^{-1} \bigg[
    -v_{n,0} \partial_z (\rho(\pf^\inn_0) \vv^\inn_0)
    + \partial_z (\rho(\pf^\inn_0) \vv^\inn_0 \otimes \vv^\inn_0) \n_0
    + \partial_z p^\inn_0 \n_0
\notag\\&
    - \partial_z \bigg(
        \mu(\pf^\inn_0) \bigg(
            \partial_z \vv^\inn_1 \otimes \n_0
            + (\partial_z \vv^\inn_1 \otimes \n_0)^T
            - \frac{2}{3} \partial_z \vv^\inn_1\cdot \n_0 \tens{I}
\notag\\&
            + \nabla_\Interface \vv^\inn_0
            + (\nabla_\Interface \vv^\inn_0)^T
            - \frac{2}{3} \nabla_\Interface \cdot \vv^\inn_0 \tens{I}
    \bigg)
\notag\\&
         + (M-1)\pf^\inn_1
        \bigg(
            \partial_z \vv^\inn_0 \otimes \n_0
            + (\partial_z \vv^\inn_0 \otimes \n_0)^T
            - \frac{2}{3} \partial_z \vv^\inn_0\cdot \n_0 \tens{I}
    \bigg) \bigg) \n_0
\notag\\&
    + \nabla_\Interface \cdot \bigg( \mu(\pf^\inn_0) \bigg(
            \partial_z \vv^\inn_0 \otimes \n_0
            + (\partial_z \vv^\inn_0 \otimes \n_0)^T
            - \frac{2}{3} \partial_z \vv^\inn_0 \cdot \n_0 \tens{I}
            \bigg) \bigg)
\notag\\&
    - \frac{3 \tension \refL}{2} \Big[
        \partial_z \Big(
            \partial_z \pf^\inn_0 \n_0
            \otimes
            (\partial_z \pf^\inn_1 \n_0 + \nabla_\Interface \pf^\inn_0 )
          +
            (\partial_z \pf^\inn_1 \n_0 + \nabla_\Interface \pf^\inn_0 )
            \otimes
            \partial_z \pf^\inn_0 \n_0
          \Big) \n_0
\notag\\&
        + \nabla_\Interface \cdot \Big(
            \partial_z \pf^\inn_0 \n_0 \otimes
            \partial_z \pf^\inn_0 \n_0
        \Big)
    + \mathcal{O}(1) \Big]
\bigg]
+ \mathcal{O}(1)
    ~.
\end{alignat}
\end{subequations}

Using the definition of the outer product, the leading order terms at $\Ox{-2}$ can be written as
\begin{equation}
\label{eq:sil_momentum_-2_rewrite}
0 =
-\frac{3 \tension \refL}{2} \partial_z (\partial_z \pf^\inn_0)^2 \n_0
-\partial_z ( \mu(\pf^\inn_0) \partial_z \vv^\inn_0)
-\partial_z \Big( \mu(\pf^\inn_0) ( (\partial_z \vv^\inn_0 \cdot \n_0)
    -\frac{2}{3} (\partial_z \vv^\inn_0 \cdot \n_0) \tens{I} )\Big) \n_0
\end{equation}
Let $\ve{t}_0$ the tangential vector, with $\ve{t}_0 \cdot \n_0 = 0$. Taking the dot product of
the above \eqref{eq:sil_momentum_-2_rewrite} with $\ve{t}_0$ yields
\begin{equation}
0 = \partial_z (\mu(\pf^\inn_0) \partial_z \vv^\inn_0) \cdot \ve{t}_0 ~.
\end{equation}
Integrating with respect to $z$ and using matching conditions, we obtain
\begin{equation}
\mu(\pf^\inn_0) \partial_z \vv^\inn_0 \cdot \ve{t}_0 = \text{\emph{const.}} = 0 ~.
\end{equation}
Since $\mu \neq 0$ we have $\partial_z \vv^\inn_0 \cdot \ve{t}_0 = 0$, corresponding to the
sharp-interface boundary condition and
together with the assumption $\partial_z \vv^\inn_0 \cdot \n_0 = 0$ the velocity $\vv^\inn_0$ is
independent of $z$.

Using also the definition of the outer product as well as $\nabla_\Interface \psi \cdot \n_0 =0$
, $\nabla_\Interface \cdot \n_0 = \kappa_0$ and $\nabla_\Interface \pf^\inn_0 = 0$, the second order
terms simplify to
\begin{alignat}{1}
0 =&
-v_{n,0} \partial_z (\rho(\pf^\inn_0) \vv^\inn_0)
    + \partial_z (\rho(\pf^\inn_0) \vv^\inn_0 \otimes \vv^\inn_0) \n_0
    + \partial_z p^\inn_0 \n_0
\notag\\&
    - \partial_z \bigg(
        \mu(\pf^\inn_0) \bigg(
            \partial_z \vv^\inn_1 \otimes \n_0
            + (\partial_z \vv^\inn_1 \otimes \n_0)^T
            - \frac{2}{3} \partial_z \vv^\inn_1\cdot \n_0 \tens{I}
            + (\nabla_\Interface \vv^\inn_0)^T
            - \frac{2}{3} \nabla_\Interface \cdot \vv^\inn_0 \tens{I}
    \bigg)
        \bigg) \n_0
\notag\\&
    - \frac{3 \tension \refL}{2} \Big[
        \partial_z \Big(
            2 \partial_z \pf^\inn_0 \partial_z \pf^\inn_1 \n_0
          \Big)
           + (\partial_z \pf^\inn_0)^2 \kappa_0 \n_0
        \Big)
     \Big]
    ~.
\end{alignat}
Integrating with respect to $z$ and applying the matching conditions yields
\begin{alignat}{1}
0 =& -v_{n,0} [\rho(\pf^\out_0) \vv^\out_0]
    + [\rho(\pf^\out_0) \vv^\out_0 \otimes \vv^\out_0] \n_0
    + [p^\out_0] \n_0
\notag\\&
    - [\mu(\pf^\out_0) ( \nabla \vv^\out_0 + (\nabla \vv^\out_0)^T - \frac{2}{3} \nabla \cdot
            \vv^\out_0 \tens{I} )\n_0 ]
    - \frac{3 \tension \refL}{2} \kappa_0 \n_0
    \int_{z=-\infty}^{\infty} (\partial_z \pf^\inn_0)^2 \dd z
~,
\end{alignat}
where $[\psi] = \psi_1 - \psi_2$ denotes the jump of $\psi$ across the interface and we used the
boundedness of
$\nabla_\Interface \vv^\inn_0$ as well as $\partial_z \pf^\inn_1$.
As the integral evaluates to $2/(3 \refL)$ and the phase velocities are equal to $v_{n,0}$ at the
interface, we recover the boundary condition \eqref{eq:si_interface_jump}.

\section{Non-dimensionalization}\label{app:non-dim}
As described in Section~\ref{sec:non-dim}, in addition to the two length scales $\ell$ and $L$ with
$\e = \ell/L$, we define reference values with dimensions
\begin{align*}
[\hat{L}] &= \mathrm{m}
&
[\hat{\ell}] &= \mathrm{m}
&
[\hat\xi] &= \mathrm{m}
&
[\hat{t}] &= \mathrm{s}
&
[\hat{v}] &= \mathrm{\frac{m}{s}}
\\
[\hat\rho] &= \mathrm{\frac{kg}{m^3}}
&
[\hat\mu] &= \mathrm{\frac{kg}{m \cdot s}}
&
[\hat{p}] &= \mathrm{\frac{kg}{ms^2}}
&
[\hat\slip] &= \mathrm{m}
&
[\hat\mob] &= \mathrm{\frac{m}{s}}
~,
\end{align*}
and let
\begin{align*}
\hat{L} &= L
&
\hat{\ell} &= \ell
&
\hat{\xi} &= \hat{\ell}
&
\hat{t} &= \frac{L}{\hat{v}}
&
\hat\rho &= \rho_2
&
\hat\mu &= \mu_2
&
\hat\slip &= \ell
&
\hat\mob &= \mob
~.
\end{align*}
This defines non-dimensionalized variables
\begin{align*}
\bar\xi &= \frac{\xi}{\hat\ell}
&
\bar\vv &= \frac{\vv}{\hat{v}}
&
\bar\rho &= \frac{\rho}{\hat\rho}
&
\bar\mu &= \frac{\mu}{\hat\mu}
&
\bar{p} &= \frac{p}{\hat{p}}
&
\bar{t} &= \frac{t}{\hat{t}}
&
\bar\slip &= \frac{\slip}{\hat\slip}
~,
\end{align*}
and dimensionless numbers
\begin{align*}
\Rey &= \frac{\hat\rho \hat{v} \hat{L}}{\hat\mu}
&
\Ca &= \frac{\hat{v} \hat\mu}{\gamma}
&
\Eu &= \frac{\hat{p}}{\hat\rho \hat{v}^2}
&
\Fr &= \frac{\hat{v}}{\sqrt{g \hat{L}}}
&
S &= \frac{\hat\mob}{\hat{v}}
~.
\end{align*}

Inserting into the mass \eqref{eq:diff_mass} and momentum \eqref{eq:diff_momentum} balances yields
\begin{subequations}
\begin{align}
\frac{\hat\rho}{\hat{t}}\frac{\partial \bar\rho}{\partial \bar{t}}
+ \frac{\hat\rho \hat{v}}{\hat{L}} \bar\nabla \cdot (\bar\rho \bar\vv) =& 0
~, \\
\frac{\hat\rho \hat{v}}{\hat{t}} \frac{\partial}{\partial \bar{t}} (\bar\rho \bar\vv)
+ \frac{\hat\rho \hat{v}^2}{\hat{L}} \bar\nabla \cdot (\bar\rho \bar\vv \otimes \bar\vv)
    =&
- \frac{\hat{p}}{\hat{L}} \bar\nabla \bar{p}
+ \frac{\hat\mu \hat{v}}{\hat{L}^2}
 \bar\nabla \cdot \left( \bar\mu \left( \bar\nabla \bar\vv + (\bar\nabla \bar\vv)^T
             - \frac{2}{3} (\bar\nabla \cdot \bar\vv)\mathbf{I} \right) \right)
\\ \notag
& + \hat\rho g \bar\rho (-\mathbf{z})
- \frac{\gamma\hat\xi}{\hat{L}^3} \frac{3\bar\xi}{2}
    \bar\nabla \cdot (\bar\nabla \bar\pf \otimes \bar\nabla \bar\pf)
~,
\end{align}
\end{subequations}
and after multiplication with $\hat{t}/\hat{\rho}$ and $\hat{L}^2/(\hat{\mu}\hat{v})$ respectively
\begin{subequations}
\begin{align}
\frac{\partial \bar{\rho}}{\partial \bar{t}}
+ \frac{\hat{t}\hat{v}}{\hat{L}} \bar{\nabla} \cdot (\bar{\rho} \bar{\vv}) =& 0
~, \\
\frac{\hat\rho\hat{L}\hat{v}}{\hat\mu}\frac{\hat{L}}{\hat{t}\hat{v}}
\frac{\partial}{\partial \bar{t}} (\bar{\rho} \bar{\vv})
+ \frac{\hat\rho\hat{L}\hat{v}}{\hat\mu}
\bar{\nabla} \cdot (\bar{\rho} \bar{\vv} \otimes \bar{\vv})
    =&
- \frac{\hat{p}\hat{L}}{\hat\mu\hat{v}} \bar{\nabla} \bar{p}
 \bar{\nabla} \cdot \left( \bar{\mu} \left( \bar{\nabla} \bar{\vv} + (\bar{\nabla} \bar{\vv})^T
             - \frac{2}{3} (\bar{\nabla} \cdot \bar{\vv})\mathbf{I} \right) \right)
\\ \notag
& - \frac{\hat\rho\hat{L}\hat{v}}{\hat\mu} \frac{g\hat{L}}{\hat{v}^2}  \bar{\rho} \mathbf{z}
- \frac{\gamma}{\hat\mu\hat{v}} \e \frac{3\xi}{2}
    \bar{\nabla} \cdot (\bar{\nabla} \bar{\pf} \otimes \bar{\nabla} \bar{\pf})
~.
\end{align}
\end{subequations}
Using the non-dimensional numbers and dropping the overline notation this is written as
\begin{subequations}
\begin{align}
\frac{\partial \rho}{\partial t} + \nabla \cdot (\rho \vv) =& 0
~, \\
\frac{\partial}{\partial t} (\rho \vv)
+ \nabla \cdot (\rho \vv \otimes \vv)
    =&
- \Eu \nabla p
+ \frac{1}{\Rey} \nabla \cdot \left( \mu \left( \nabla \vv + (\nabla \vv)^T
             - \frac{2}{3} (\nabla \cdot \vv)\mathbf{I} \right) \right)
\\ \notag
& - \frac{1}{\Fr^2} \rho \mathbf{z}
- \frac{\e}{\Ca} \frac{3\xi}{2}
    \nabla \cdot (\nabla \pf \otimes \nabla \pf)
~.
\end{align}
\end{subequations}
For the slip boundary condition \eqref{eq:diff_bc_slip} inserting the non-dimensionalized variables
yields
\begin{equation}
\hat{v} \bar{\vv} = \frac{\hat{\lambda}\hat{v}}{\hat{L}} \bar{\lambda} (\bar{\partial}_\n \bar{\vv}_\ve{t})\ve{t}
~, \text{ on } \InterfaceCell
~.
\end{equation}
With $\hat{\lambda} = \hat{\ell}$ and dropping the overline notation they simplify to
\begin{equation}
\vv = \e \lambda (\partial_\n \vv_\ve{t})\ve{t}
~, \text{ on } \InterfaceCell
~.
\end{equation}

\subsection{Phase field}
The phase-field equation \eqref{eq:pore_phasefield}
\begin{equation}
\xi^2 \frac{\partial \pf}{\partial t} - \mob \xi^3 \nabla^2 \pf + \xi^2 \nabla \cdot (\vv \pf) =
-\mob \xi P'(\pf)
\end{equation}
can be rewritten with non-dimensionalized variables as
\begin{equation}
\frac{\hat\xi^2}{\hat{t}} \bar\xi^2 \partial_{\bar{t}} \bar\pf
- \frac{\hat\mob\hat\xi^3}{\hat{L}^2} \bar\xi^3 \bar\nabla^2 \bar\pf
+ \frac{\hat\xi^2\hat{v}}{\hat{L}} \bar\xi^2 \bar\nabla \cdot (\bar\vv \bar\pf)
= -\hat\mob \hat\xi \bar\xi P'(\bar\pf)
~.
\end{equation}
Dividing by $\hat\xi^2/\hat{t}$, using the dimensionless number $S$ and dropping the overline notation the equation becomes
\begin{equation}
\xi^2 \partial_t \pf - S \e^1 \xi^2 \nabla^2 \pf + \xi^2 \nabla \cdot (\vv \pf)
 = - S \e^{-1} P'(\pf)
~.
\end{equation}
For the contact-angle boundary condition we obtain
\begin{equation}
\frac{1}{\hat{L}} \bar{\partial}_\n \bar{\pf} =
- \frac{1}{\hat{\xi}} \cos(\contactAngle) \bar{\xi}^{-1} \sqrt{2P(\bar{\pf})}
~, \text{ on } \InterfaceCell
~.
\end{equation}
Multiplying with $L$, using $\hat{\xi} = \hat{\ell} = \e \hat{L}$ and dropping the overline notation this yields
\begin{equation}
\partial_\n \pf =
- \e^{-1} \cos(\contactAngle) \xi^{-1} \sqrt{2P(\pf)}
~, \text{ on } \InterfaceCell
~.
\end{equation}

\section{Homogenization}\label{app:homogenization}
We introduce the micro-scale coordinate $\y = \e^{-1}\x$ with $\e = \ell/L$ the scale separation and
assume asymptotic expansions for $\psi \in \{p, \vv, \pf\}$
\begin{equation}
\psi(t,\x) = \sum_{k=0}^\infty \e^k \psi_k \left(t, \x, \frac{\x}{\e}\right)
\end{equation}
with each $\psi_k$ $y$-periodic on the reference cell $Y$.
Rewriting spatial derivatives according to
\begin{equation}
\nabla\psi = \nabla_\x \sum_{k=0}^\infty \e^k \psi_k (t, \x, \y)
    + \frac{1}{\e} \nabla_\y \sum_{k=0}^\infty \e^k \psi_k (t, \x, \y)
\end{equation}
and choosing non-dimensional numbers $\Ca\in\Oe{0}$, $\Rey\in\Oe{0}$, $\Eu\in\Oe{-2}$ and
$\Fr\in\Oe{0}$ in addition
to $S \in\Oe{0}$ and inserting the asymptotic expansions into the non-dimensionalized equations
yields the following leading order terms.
For the mass balance \eqref{eq:non-dim_mass} we obtain
\begin{equation}
\begin{aligned}
0 &= \frac{\partial \rho}{\partial t} + \nabla \cdot (\rho \vv)
\\&= \frac{\partial}{\partial t} (\rho_0 + \Oe{1}) + (\e^{-1} \nabla_\y + \nabla_\x)
    \cdot ((\rho_0 + \e \rho_1 + \Oe{2}) (\vv_0 + \e \vv_1 + \Oe{2}))
\\&= \bigg(\frac{\partial \rho_0}{\partial t} + \Oe{1}\bigg)
    + (\e^{-1} \nabla_\y + \nabla_\x) \cdot
    (\rho_0 \vv_0 + \e (\rho_0 \vv_1 + \rho_1 \vv_0) + \Oe{2})
\\&= \bigg[\frac{\partial \rho_0}{\partial t}\bigg] + \Oe{1}
    + \e^{-1} \big[\nabla_\y \cdot (\rho_0 \vv_0) \big]
    + \big[ \nabla_\y \cdot (\rho_0 \vv_1 + \rho_1 \vv_0)
    + \nabla_\x \cdot (\rho_0 \vv_0) \big]
    + \Oe{1}
\\&= \e^{-1} \big[ \nabla_\y \cdot (\rho_0 \vv_0) \big]
    +\e^{0} \bigg[ \frac{\partial \rho_0}{\partial t}
    + \nabla_\y \cdot (\rho_0 \vv_1 + \rho_1 \vv_0)
    + \nabla_\x \cdot (\rho_0 \vv_0) \bigg]
    + \Oe{1}
~.
\end{aligned}
\end{equation}
We note that due to the dependence of $\rho$ on $\pf$ we have
\begin{equation}
\rho_0 = \rho(\pf_0) ~, \qquad \rho_1 = \pf_1 (R - 1) ~,
\end{equation}
and the equation can be written as
\begin{equation}
0 =  \e^{-1} \big[ \nabla_\y \cdot (\rho(\pf_0) \vv_0) \big]
    +\e^{0} \bigg[ \frac{\partial \rho(\pf_0)}{\partial t}
    + \nabla_\y \cdot (\rho(\pf_0) \vv_1 + \rho_1 \vv_0)
    + \nabla_\x \cdot (\rho(\pf_0) \vv_0) \bigg]
    + \Oe{1}
~.
\end{equation}
Analogously for the viscosity $\mu$ we have $\mu_0 = \mu(\pf_0)$ and $\mu_1 = \pf_1 (M-1)$.
Denoting $\overline{\Eu} := \e^2 \Eu \in \Oe{0}$, the momentum equation \eqref{eq:non-dim_momentum} yields
the following terms.
\begin{subequations}
\begin{alignat}{1}
0 =& -\frac{\partial}{\partial t} (\rho \vv)
     - \nabla \cdot (\rho \vv \otimes \vv)
     - \Eu \nabla p
     + \frac{1}{\Rey} \nabla \cdot \left( \mu \left( \nabla \vv + (\nabla \vv)^T
     - \frac{2}{3} (\nabla \cdot \vv)\mathbf{I} \right) \right)
\notag\\&
     - \frac{1}{\Fr^2} \rho \mathbf{z}
     - \frac{\e}{\Ca} \frac{3\xi}{2} \nabla \cdot (\nabla \pf \otimes \nabla \pf)
\\=& -\frac{\partial}{\partial t} ((\rho_0 + \Oe{1})(\vv_0 + \Oe{1}))
\notag\\&
     - (\e^{-1} \nabla_\y + \nabla_\x) \cdot
        ((\rho_0 + \Oe{1}) (\vv_0 + \Oe{1}) \otimes (\vv_0 + \Oe{1}))
\notag\\&
     - \e^{-2} \overline{\Eu} (\e^{-1} \nabla_\y + \nabla_\x) (p_0 + \e p_1 + \Oe{2})
\notag\\&
     + \frac{1}{\Rey} (\e^{-1} \nabla_\y + \nabla_\x) \cdot
     \bigg( (\mu_0 + \Oe{1}) \bigg( (\e^{-1} \nabla_\y + \nabla_x) (\vv_0 + \Oe{1})
\notag\\&
        + ((\e^{-1} \nabla_\y + \nabla_x) (\vv_0 + \Oe{1}))^T
    - \frac{2}{3} \big((\e^{-1} \nabla_\y + \nabla_\x) \cdot (\vv_0 + \Oe{1})\big)\mathbf{I} \bigg) \bigg)
\notag\\&
    - \frac{1}{\Fr^2} (\rho_0 + \Oe{1}) \mathbf{z}
\notag\\&
    - \frac{\e}{\Ca} \frac{3\xi}{2} (\e^{-1} \nabla_\y + \nabla_\x) \cdot
    \big((\e^{-1}\nabla_y + \nabla_x) (\pf_0 + \Oe{1}) \otimes (\e^{-1}\nabla_y + \nabla_x) (\pf_0 + \Oe{1})\big)
\\=& -\Oe{0} -\Oe{-1}
     - \e^{-3} \big[\overline{\Eu} \nabla_\y p_0 \big]
     - \e^{-2} \big[ \overline{\Eu}(\nabla_\y p_1 + \nabla_\x p_0) \big] - \Oe{-1}
\notag\\&
     + \e^{-2} \bigg[ \frac{1}{\Rey} \nabla_y \cdot \left( \mu_0 \left( \nabla_\y \vv_0 + (\nabla_\y \vv_0)^T
                 -\frac{2}{3} (\nabla_\y \cdot \vv_0)\mathbf{I} \right)\right) \bigg]
     + \Oe{-1}
     - \Oe{0}
\notag\\&
     - \e^{-2} \bigg[ \frac{1}{\Ca} \frac{3\xi}{2} \nabla_\y \cdot
        (\nabla_\y \pf_0 \otimes \nabla_\y \pf_0) \bigg]
     - \Oe{-1}
\\=& \,\e^{-3} \big[-\overline{\Eu} \nabla_\y p_0 \big]
\notag\\& +\e^{-2} \bigg[
        \overline{\Eu}(-\nabla_\y p_1 - \nabla_\x p_0)
        +\frac{1}{\Rey} \nabla_y \cdot \left( \mu(\pf_0) \left( \nabla_\y \vv_0 + (\nabla_\y \vv_0)^T
                 -\frac{2}{3} (\nabla_\y \cdot \vv_0)\mathbf{I} \right)\right)
\notag\\&
        -\frac{1}{\Ca} \frac{3\xi}{2} \nabla_\y \cdot
        (\nabla_\y \pf_0 \otimes \nabla_\y \pf_0)
    \bigg]
\end{alignat}
\end{subequations}
Using the polynomial structure of $P'$ we obtain the following expansion from the phase-field
equation \eqref{eq:non-dim_pf}.
\begin{subequations}
\begin{alignat}{1}
0 =& \xi^2 \frac{\partial \pf}{\partial t} - S \e^1 \xi^2 \nabla^2 \pf + \xi^2 \nabla \cdot (\vv \pf) + S \e^{-1} P'(\pf)
\\=& \xi^2 \frac{\partial}{\partial t} (\pf_0 + \Oe{1})
    - \e^1 S\xi^2 (\e^{-1}\nabla_\y + \nabla_\x) \cdot
     (\e^{-1}\nabla_\y + \nabla_\x) (\pf_0 + \e \pf_1 + \Oe{2})
\notag\\&
    + \xi^2 (\e^{-1}\nabla_\y + \nabla_\x) \cdot
     ((\vv_0 + \e\vv_1 + \Oe{2}) (\pf_0 + \e\pf_1 + \Oe{2}))
\notag\\&
    + \e^{-1} S (P'(\pf_0) + \e \pf_1 P''(\pf_0) + \Oe{2})
\\=& \bigg[ \xi^2 \frac{\partial \pf_0}{\partial t} \bigg] + \Oe{1}
\notag\\&
    - \e^{-1} \big[ S \xi^2 \nabla_\y^2 \pf_0 \big]
    - \big[ S \xi^2 (\nabla_\y \cdot (\nabla_\x \pf_0) + \nabla_\x \cdot (\nabla_\y \pf_0) +
            \nabla_\y \cdot (\nabla_\y \pf_1)) \big] -\Oe{1}
\notag\\&
    + \e^{-1} \big[ \xi^2 \nabla_\y \cdot (\vv_0 \pf_0) \big]
    + \e^{0} \big[ \xi^2 \nabla_\y \cdot (\vv_0 \pf_1 + \vv_1 \pf_0)
    + \xi^2\nabla_\x \cdot (\vv_0 \pf_0) \big]
\notag\\&
    + \e^{-1} \big[ S P'(\pf_0) \big]
    + \big[ S \pf_1 P''(\pf_0) \big]
\\=& \e^{-1} \big[ -S\xi^2 \nabla_\y^2 \pf_0
    + \xi^2 \nabla_\y \cdot (\vv_0 \pf_0) + S P'(\pf_0) \big]
\notag\\&
    + \e^{0} \bigg[ \xi^2 \frac{\partial \pf_0}{\partial t}
    + \xi^2 \nabla_\y \cdot (\vv_0 \pf_1 + \vv_1 \pf_0) + \xi^2 \nabla_\x \cdot (\vv_0 \pf_0)
\notag\\&
    - S \xi^2 ( \nabla_\y \cdot (\nabla_\x \pf_0) + \nabla_\x \cdot (\nabla_\y \pf_0) +
            \nabla_\y \cdot (\nabla_\y \pf_1) )
    + S \pf_1 P''(\pf_0) \bigg]
\end{alignat}
\end{subequations}

Inserting the asymptotic expansions into the boundary conditions yields

\begin{subequations}
\begin{alignat}{2}
0 &= \vv - \e \lambda (\partial_\n \vv_\ve{t})\ve{t}
\notag\\&= (\vv_0 + \Oe{1})
    - \e \lambda ((\e^{-1}\nabla_y + \nabla_\x)(\vv_0 + \Oe{1}) \n )
\notag\\&= \e^{0} \big[ \vv_0 - \lambda \nabla_y \vv_0 \n \big] + \Oe{1}
~, && \text{ on } \InterfaceCell
~, \\
0 &= \partial_\n \pf + \e^{-1} \cos(\contactAngle) \xi^{-1} \sqrt{2P(\pf)}
\notag\\&= (\e^{-1} \nabla_\y + \nabla_\x) (\pf_0 + \Oe{1}) \n
    + \e^{-1} \cos(\contactAngle) \xi^{-1} (\sqrt{2P(\pf_0)} + \Oe{1})
\notag\\&= \e^{-1} \big[ \nabla_y \pf_0 \n + \cos(\contactAngle) \xi^{-1} \sqrt{2P(\pf_0)} \big]
        + \Oe{1}
~, && \text{ on } \InterfaceCell
~.
\end{alignat}
\end{subequations}

\bibliography{preprint_kelm_20230630}

\end{document}